\documentclass[11pt,a4paper]{article} 
\usepackage{graphicx,amsmath,bm}
\usepackage{amssymb,amscd}
\usepackage{caption2,color}

\usepackage[dvipdfm,
           pdfstartview=FitH,
           CJKbookmarks=true,
           bookmarksnumbered=true,
           bookmarksopen=true,
          colorlinks=true,
          colorlinks=cyan,
           pdfborder=001,
           citecolor=blue%
           ]{hyperref}

\makeatletter

\@addtoreset{equation}{section} \makeatother

\newtheorem{lemma}{Lemma}[section]

\newtheorem{rem}{Remark}[section]
\newtheorem{theorem}{Theorem}[section]
\newtheorem{corollary}{Corollary}[section]
\newcommand{\dx}{\,\mathrm{d}x}

\newcommand{\ds}{\,\mathrm{d}s}
\newcommand{\n}{\nabla}
\newcommand{\p}{\partial}
\newcommand{\norm}[1]{\lVert#1\rVert}
\newcommand{\seminorm}[1]{\lvert#1\rvert}
\newcommand{\rn}{\ensuremath{\mathbb{R}^N}}
\DeclareMathAlphabet{\mathsfsl}{OT1}{cmss}{m}{sl}
\newcommand{\tensor}[1]{\mathsfsl{#1}}
\renewcommand{\vec}[1]{\mbox{\boldmath$#1$}}
\newcommand{\oo}{\ensuremath{\Omega}}
\newcommand{\G}{\Gamma}
\newcommand{\ve}{\varepsilon}

\newcommand{\diff}{\,\mathrm{d}}

\newcommand{\mdiv}{\,\mathrm{div}\,}

\newcommand{\md}{\mathrm{D}}
\newcommand{\mg}{\mathcal{G}}
\newcommand{\mt}{\mathcal{T}}

\newcommand{\vv}{\vec{V}}

\newcommand{\vvv}{\vec{v}}

\newcommand{\vy}{\vec{y}}

\newcommand{\vu}{\vec{u}}

\newcommand{\vg}{\vec{g}}

\newcommand{\vphi}{\vec{\varphi}}

\newcommand{\vn}{\vec{n}}
\begin{document}

\title{Optimal Shape Design for the Viscous Incompressible Flow\footnote{This work was
supported by the National Natural Science Fund of China under grant
number 10671153.}}

\author{Zhiming Gao\thanks{Corresponding author.  School of Science, Xi'an Jiaotong University, P.O.Box 1844, Xi'an, Shaanxi, P.R.China, 710049. E--mail:dtgaozm@gmail.com.}\qquad
and \qquad Yichen Ma\footnote{School of Science, Xi'an Jiaotong
University, Shaanxi, P.R.China, 710049.
E-mail:\,ycma@mail.xjtu.edu.cn.}
 }

 \maketitle
\noindent{{\textbf{Abstract.\;}} This paper is concerned with a
numerical simulation of shape optimization in a two-dimensional
viscous incompressible flow governed by Navier--Stokes equations
with mixed boundary conditions containing the pressure. The
minimization problem of total dissipated energy was established in
the fluid domain. We derive the structures of shape gradient of the
cost functional by using the differentiability of a minimax
formulation involving a Lagrange functional with a function space
parametrization technique.
 Finally a gradient type
algorithm is
effectively used for our problem.   \\[8pt]
{{\textbf{Keywords.\;}}
shape optimization; minimax principle; gradient algorithm; Navier--Stokes equations.\\[8pt]
{{\textbf{AMS (2000) subject classifications.\;}}35B37, 35Q30,49K35.

\section{Introduction}
The problem of finding the optimal design of a system for the
viscous incompressible flow arises in many design problems in
aerospace, automotive, hydraulic, ocean, structural, and wind
engineering. In practise, engineers are interested in reducing the
drag force in the wing of a plane or vehicle or in reducing the
dissipation in channels, hydraulic values, etc.

The optimal shape design of a body subjected to the minimum
dissipate viscous energy has been a challenging task for a long
time, and it has been investigated by several authors. For instance,
O.Pironneau in \cite{Piron,piron01} computes  the derivative of the
cost functional using normal variation approach; F.Murat and J.Simon
in \cite{murat74} use the formal calculus to deduce an expression
for the derivative; J.A.Bello \emph{et.al.} in
\cite{bello92a,bello92b,bello97} considered this problem
theoretically in the case of Navier--Stokes flow by the formal
calculus.

The previous references concern Dirichlet boundary conditions for
the velocity. However, the velocity--pressure type boundary
conditions must be introduced and it seems more realistic in many
industrial applications, such as shape optimization of
Aorto--Coronaric bypass anastomoses in biomedical engineering.
Recently, H.Yagi and M. Kawahara in \cite{yagi05} study the optimal
shape design for Navier--Stokes flow with boundary conditions
containing the pressure using a discretize--then--differentiate
approach. However, its proposed algorithm converges slowly.
E.Katamine \emph{et.al.}\cite{azeg05} use the
differentiate--then--discretize approach with the formulae of
material derivative to study such problem involving
velocity--pressure type boundary conditions with Reynolds number up
to 100.

In the present paper, we will use the so-called function space
parametrization technique which was advocated by M.C.Delfour and
J.-P.Zol\'{e}sio to solving poisson equation with Dirichlet and
Nuemann condition (see \cite{delfour}). In our paper
\cite{gao0601,gao06,gao06ns}, we apply them to solve a Robin problem
and shape reconstruction problems for Stokes and Navier--Stokes flow
with Dirichlet boundary condition only involving the velocity,
respectively.

 However,
in this paper we extend them to study the energy minimization
problem for Navier--Stokes flow with velocity--pressure boundary
conditions in despite of its lack of rigorous mathematical
justification in case where the Lagrange formulation is not convex.
We shall show how this theorem allows, at least formally to bypass
the study of material derivative and obtain the expression of shape
gradient for the given cost functional. Finally we will introduce an
efficient numerical algorithm for the solution of such minimization
problems and the numerical examples show that our proposed algorithm
converges very fast.

This paper is organized as follows. In section 2, we briefly recall
the velocity method which is used for the characterization of the
deformation of the shape of the domain and give the description of
the shape minimization problem for the Navier--Stokes flow with
mixed type boundary conditions.

Section 3 is devoted to the computation of the shape gradient of the
Lagrangian functional due to a minimax principle concerning the
differentiability of the minimax formulation by function space
parametrization technique.

Finally in section 4, we give its finite element approximation and
propose a gradient type algorithm with some numerical examples to
show that our theory could be very useful for the practical purpose
and the proposed algorithm is efficient.
\section{Preliminaries and statement of the problem}\label{sec2}
\subsection{Elements of the velocity method}
To our little knowledge, there are about three types of techniques
to perform the domain deformation: J.Hadamard \cite{ha07}'s normal
variation method, the {perturbation of the identity} method by
J.Simon \cite{si80} and the {velocity method} (see J.Cea\cite{ce81}
and J.-P.Zolesio\cite{delfour,zo79}). We will use the velocity
method which contains the others. In that purpose, we choose an open
set $D$ in $\rn$ with the boundary $\p D$ piecewise $C^k$, and a
velocity space $\vec V\in \mathrm{E}^k
:=C([0,\varepsilon];[\mathcal{D}^k(\bar{D})]^N)$, where
$\varepsilon$ is a small positive real number and
$[\mathcal{D}^k(\bar{D})]^N$ denotes the space of all $k-$times
continuous differentiable functions with compact support contained
in $D$. The velocity field
$$\vec V(t)(x)=\vec V(t,x), \qquad x\in D,\quad t\geq 0$$
belongs to $[\mathcal{D}^k(\bar{D})]^N$ for each $t$. It can
generate transformations $T_t(\vec V)X=x(t,X)$ through the following
dynamical system
\begin{equation*}
  \frac{\diff x}{\diff t}(t,X)=\vec V(t,x(t)),\qquad  x(0,X)=X
\end{equation*}
with the initial value $X$ given. We denote the "transformed domain"
$T_t(\vec V)(\oo)$ by $\oo_t(\vec V)$ at $t\geq 0$, and also set its
boundary $\Gamma_t:=T_t(\Gamma)$.

There exists an interval $I=[0,\delta)$, $0<\delta\leq\varepsilon,$
and a one-to-one map $T_t$ from $\bar{D}$ onto $\bar{D}$ such that
\begin{itemize}
    \item [(i)] $T_0=\mathrm{I};$
    \item [(ii)] $(t,x)\mapsto T_t(x)$ belongs to $C^1(I;C^k(D))$ with $T_t(\p D)=\p D$;
    \item[(iii)]$(t,x)\mapsto T_t^{-1}(x)$ belongs to $C(I;C^k(D))$.
\end{itemize}
Such transformation are well studied in \cite{delfour}.

 Furthermore, for sufficiently small $t>0,$ the Jacobian $J_t$ is
 strictly positive:
 \begin{equation*}\label{jacobian}
   J_t(x):=\det\seminorm{\md T_t(x)}=\det\md T_t(x)>0,
 \end{equation*}
where $\md T_t(x)$ denotes the Jacobian matrix of the transformation
$T_t$ evaluated at a point $x\in D$ associated with the velocity
field $\vec V$. We will also use the following notation: $\md
T_t^{-1}(x)$ is the inverse of the matrix $\md T_t(x)$ , ${}^*\md
T_t^{-1}(x)$ is the transpose of the matrix $\md T_t^{-1}(x)$. These
quantities also satisfy the following lemma.
\begin{lemma}\label{lem:a}
    For any $\vec V\in E^k$, $\md T_t$ and $J_t$ are invertible. Moreover, $\md T_t$, $\md T_t^{-1}$
    are in $C^1([0,\varepsilon];[C^{k-1}(\bar{D})]^{N\times N})$,
    and $J_t$, $J_t^{-1}$ are in
    $C^1([0,\varepsilon];C^{k-1}(\bar{D}))$.
\end{lemma}

Now let $J(\oo)$ be a real valued functional associated with any
regular domain $\oo$, we say that the functional $J(\oo)$ has a {\bf
Eulerian derivative} at
$\oo$ in the direction $\vec V$ if the limit\\[6pt]
\begin{equation*}
\lim_{t\searrow 0}\frac{1}{t}\left[J(\oo_t)-J(\oo)\right]:=\diff
J(\oo;\vec V)
\end{equation*}
exists.

 Furthermore, if the map $\vec V\mapsto\diff
J(\oo;\vec V):\;\mathrm{E}^k\rightarrow\mathbb{R}$ is linear and
continuous, we say that $J$ is {\bf shape differentiable} at $\oo$.
In the distributional sense we have
\begin{equation}\label{pri:shaped}
    \diff J(\oo;\vec V)=\langle \n J,\vec V\rangle_{(\mathcal{D}^k(\bar{D})^N)'\times \mathcal{D}^k(\bar{D})^N}.
\end{equation}
 When $J$ has a Eulerian derivative, we say that $\n J$ is the {\bf shape gradient} of $J$
at $\oo$.

Before closing this subsection, we introduce the following
functional spaces which will be used in this paper:
\begin{eqnarray*}
H(\oo)&:=&\{\vu\in  H^1(\oo)^N: \;\mdiv\vu=0\mbox{ in }\oo,\;\vu=0\mbox{ on }\G_w\cup\G_s,\;\vu=\vg\mbox{ on }\G_u\},\\
V_g(\oo)&:=&\{\vu\in  H^2(\oo)^N: \;\vu=0\mbox{ on }\G_w\cup\G_s,\;\vu=\vg\mbox{ on }\G_u\},\\
V_0(\oo)&:=&\{\vu\in H^2(\oo)^N: \;\vu=0\mbox{ on }\G_w\cup
\G_u\cup \G_s\},\\
 Q(\oo)&:=&\left\{p\in H^1(\oo):\;\int_\oo p\dx=0 \;(\mbox{ if
meas}(\G_d)=0)\right\}.
 \end{eqnarray*}
\subsection{Formulation of the flow optimization problem}
Let $\oo$ be a region of $\mathbb{R}^2$ and we denote by $\G$ the
boundary of $\oo$. We suppose that $\oo$ is filled with a Newtonian
incompressible viscous fluid of the kinematic viscosity $\nu$. The
flow of such a fluid is modeled by the following system of
Navier--Stokes equations,
\begin{eqnarray}
    \label{nsp:a}   & -\mdiv\tensor{\sigma}(\vy,p)+\md\vy\cdot\vy=0\; &\qquad \mbox{in }\oo,\\
     \label{nsp:b}   &\mdiv\vy=0 & \qquad\mbox{in }\oo,
\end{eqnarray}
where $\vy$ denotes the velocity field, $p$ the pressure, and
$\sigma(\vy,p)$ the stress tensor defined by
$\sigma(\vy,p):=-p\mathrm{I}+2\nu\ve(\vy)$ with the rate of
deformation tensor $\varepsilon(\vy):=(\md\vy+{ }^*\md\vy)/2,$ where
${}^*\md\vy$ denotes the transpose of the matrix $\md\vy$ and
$\mathrm{I}$ denotes the identity tensor.

 Equations \eqref{nsp:a} and \eqref{nsp:b} have to be completed by
 the following typical
boundary conditions:
\begin{eqnarray}
\label{nsp:c}  \vy=\vg &&\quad \mbox{on }\Gamma_u\\
\label{nsp:d}     \vy=0 &&\quad\mbox{on }\G_s\cup \Gamma_w \\
\label{nsp:e}     \sigma(\vy,p)\cdot\vn=\vec h &&\quad\mbox{on
}\Gamma_d
\end{eqnarray}
where $\vn$ denotes the unit vector of outward normal on
$\G=\G_u\cup\G_d\cup\G_w\cup\G_s$, $\G_u$ is the inflow boundary,
$\G_d$ the outflow boundary, $\G_w$ the boundary corresponding to
the fluid wall and $\G_s$ is the boundary which is to be optimized.
We also recall that the Reynolds number $\mathrm{Re}$ is classically
defined by $\mathrm{Re}=UL/\nu$ with $U$ a characteristic velocity
and $L$ a characteristic length.

For the existence and uniqueness of the solution of the
Navier--Stokes system \eqref{nsp:a}--\eqref{nsp:e}, we have the
following results (see \cite{temam01}).
\begin{theorem}
        \label{thm:ns}
        Suppose that $\oo$ is of class $C^1$. For the data
\begin{equation*}
\vg\in H^{3/2}(D)^N,\qquad\int_{\G_u}\vg\cdot\vn\ds=0;\qquad \vec
h\in H^{1/2}(D)^N,
\end{equation*}
 there exists at least one $\vy\in H(\oo)$ and a distribution $p\in L^2(\oo)$
 on $\oo$ such that \eqref{nsp:a}--\eqref{nsp:e} hold. Moreover,
 if $\nu$ is sufficiently large or $\vg$ and $\vec h$ sufficiently small,
 there exists a unique solution $(\vy,p)\in H(\oo)\times L^2(\oo)$ to
 the problem \eqref{nsp:a}--\eqref{nsp:e}. In addition, if $\oo$ is of class $C^2$, we
 have $(\vy, p)\in V_g(\oo)\times
Q(\oo)$.
\end{theorem}
Our goal is to optimize the shape of the domain $\oo$ which
minimizes a given cost functional depending on the fluid domain and
the state. The cost functional may represent a given objective
related to specific characteristic features of the flow (e.g., the
deviation with respect to a given target pressure, the drag, the
vorticity, $\cdots$).

 We are interested in
solving the total dissipation energy minimization problem
\begin{equation}\label{nsdrag:cost}
 \min_{\oo\in\mathcal{O}} J(\oo)=2\nu\int_{\oo}\seminorm{\varepsilon(\vy)}^2\dx,
\end{equation}
where the boundary $\G_u\cup\G_d\cup\G_w$ is fixed and an example of
the admissible set ${\mathcal{O}}$ is:
$$\mathcal{O}:=\left\{\oo\subset\rn:\; \G_u\cup\G_d\cup\G_w \mbox{ is
fixed},\;\int_\oo\dx=\mbox{constant}\right\}.$$

\begin{corollary}[\cite{piron98}]
  Let $\oo$ be of piecewise $C^1$, the minimization problem \eqref{nsdrag:cost} has at least one
  solution with given area in two dimensions.
\end{corollary}

\setcounter{figure}{0}
\section{Function space parametrization}\label{sec4}
In this section we derive the structure of the shape gradient for
the cost functional $J(\oo)$ by function space parametrization
techniques in order to bypass the usual study of material
derivative.

Let $\oo$ be of class $C^2$, the weak formulation of
\eqref{nsp:a}--\eqref{nsp:e} in mixed form is:
\begin{equation}\label{state:weak}
  \left\{
  \begin{array}
    {ll}
    &\mbox{seek } (\vy,p)\in V_g(\oo)\times Q(\oo)\mbox{ such
  that}\\[4pt]
  &   \int_{\oo}[2\nu\varepsilon(\vy):\varepsilon(\vvv)
     +\md\vy\cdot\vy\cdot\vvv-p\mdiv\vvv]\dx=\int_{\G_d}\vec
     h\cdot\vvv\ds,\;\forall \vvv\in V_0(\oo),\\[4pt]
 &\int_{\oo}\mdiv\vy q\dx=0,\;\forall q\in Q(\oo).
  \end{array}
  \right.
\end{equation}
Where in the weak form \eqref{state:weak}, we have used the
following lemma.
\begin{lemma}  \label{lem:green}
\begin{equation*}
  2\int_\oo\ve(\vy):\ve(\vvv)\dx=-\int_\oo(\Delta\vy+\n\mdiv\vy)\cdot\vvv\dx+2\int_{\p\oo}\ve(\vy)\cdot\vn\cdot\vvv\ds.
\end{equation*}
\end{lemma}
Now we introduce the following Lagrange functional associated with
\eqref{state:weak} and \eqref{nsdrag:cost}:
\begin{equation}
  G(\oo,\vy,p,\vvv,q)=J(\oo)-L(\oo,\vy,p,\vvv,q),
\end{equation}
where
\begin{equation*}
  L(\oo,\vy,p,\vvv,q)=\int_{\oo}[2\nu\varepsilon(\vy):\varepsilon(\vvv)
     +\md\vy\cdot\vy\cdot\vvv-p\mdiv\vvv]\dx-\int_{\G_d}\vec
     h\cdot\vvv\ds-\int_{\oo}\mdiv\vy q\dx.
\end{equation*}
The minimization problem \eqref{nsdrag:cost} can be put in the
following form
\begin{equation}
  \min_{\oo\in\mathcal{O}}\;\min_{(\vy,p)\in V_g(\oo)\times Q(\oo)}\;\max_{(\vvv,q)\in V_0(\oo)\times
  Q(\oo)}G(\oo,\vy,p,\vvv,q),
\end{equation}
 We can use the minimax framework to avoid the study of the
state derivative with respect to the shape of the domain. The
Karusch-Kuhn-Tucker conditions will furnish the shape gradient of
the cost functional $J(\oo)$ by using the adjoint system. Now let's
establish the first optimality condition for the problem
\begin{equation}
\min_{(\vy,p)\in V_g(\oo)\times Q(\oo)}\;\max_{(\vvv,q)\in
V_0(\oo)\times
  Q(\oo)}G(\oo,\vy,p,\vvv,q).
\end{equation}
Formally the adjoint equations are defined from the Euler--Lagrange
equations of the Lagrange functional $G$. Clearly, the variation of
$G$ with respect to $(\vvv,q)$ can recover the state system
\eqref{state:weak}. On the other hand, in order to find the adjoint
state system, we differentiate $G$ with respect to $p$ in the
direction $\delta p$,
$$  \frac{\p
  G}{\p p}(\oo,\vy,p,\vvv,q)\cdot\delta p=\int_\oo\delta p\mdiv\vvv\dx=0,$$
Taking $\delta p$ with compact support in $\oo$ gives
  \begin{equation}
     \mdiv\vvv=0.\label{adj:b}
  \end{equation}
Then we differentiate $G$ with respect to $\vy$ in the direction
$\delta\vy$ and employ Green formula,
\begin{multline*}
  \frac{\p
  G}{\p\vy}(\oo,\vy,p,\vvv,q)\cdot\delta\vy=\int_\oo(-2\nu\Delta\vy+\nu\Delta\vvv-\n
  q-{}^*\md\vy\cdot\vvv+\md\vvv\cdot\vy)\cdot\delta\vy\dx\\
-\int_{\p\oo}\sigma(\vvv,q)\cdot\vn\cdot\delta\vy\ds+4\nu\int_{\p\oo}\ve(\vy)\cdot\vn\cdot\delta\vy\ds
  -\int_{\p\oo}(\vy\cdot\vn)(\vvv\cdot\delta\vy)\ds.
\end{multline*}
Taking $\delta\vy$ with compact support in $\oo$ gives
\begin{equation}\label{adj:a}
-\nu\Delta\vvv+\n
  q+{}^*\md\vy\cdot\vvv-\md\vvv\cdot\vy=-2\nu\Delta\vy.
  \end{equation}
Then varying $\delta\vy$ on $\G_d$ gives
\begin{equation}
  \sigma(\vvv,q)\cdot\vn+(\vy\cdot\vn)\vvv-4\nu\ve(\vy)\cdot\vn=0,\qquad\mbox{on
  }\G_d.
\end{equation}
Finally we obtain the following adjoint state system
\begin{equation}
  \left\{
  \begin{array}
    {ll}
-\mdiv\sigma(\vvv,q)+{}^*\md\vy\cdot\vvv-\md\vvv\cdot\vy=-2\nu\Delta\vy
&\qquad\mbox{ in
  }\oo\\
  \mdiv\vvv=0 &\qquad\mbox{  in }\oo\\
  \sigma(\vvv,q)\cdot\vn+(\vy\cdot\vn)\vvv-4\nu\ve(\vy)\cdot\vn=0,&\qquad\mbox{ on
  }\G_d\\
  \vvv=0&\qquad\mbox{ on }\G_u\cup\G_w\cup\G_s,
  \end{array}
  \right.
\end{equation}
and its variational form
\begin{equation}\label{adj:weak}
  \left\{
  \begin{array}
    {ll}
    &\mbox{seek } (\vvv,q)\in V_0(\oo)\times Q(\oo)\mbox{ such
  that}\;\forall (\vphi,\psi)\in V_0(\oo)\times Q(\oo),\\[5pt]
  &   \int_{\oo}[2\nu\varepsilon(\vvv):\varepsilon(\vphi)
     +\md\vphi\cdot\vy\cdot\vvv+\md\vy\cdot\vphi\cdot\vvv-q\mdiv\vphi]\dx=4\nu\int_\oo\ve(\vy):\ve(\vphi)\dx,\\[4pt]
 &\int_{\oo}\mdiv\vvv \psi\dx=0.
  \end{array}
  \right.
\end{equation}
We employ the velocity method to modelize the domain deformations.
We only perturb the boundary $\G_s$ and consider the mapping
$T_t(\vv)$, the flow of the velocity field $$\vv\in
V_{\mathrm{ad}}:=\{\vv\in C^0(0,\tau;C^2(\rn)^N):\,V=0\;\mbox{in the
neighorhood of } \G_u\cup\G_w\cup\G_d\}.$$ We denote the perturbed
domain $\oo_t=T_t(\vv)(\oo)$.

Our objective in this section is to study the derivative of $j(t)$
with respect to $t$, where
\begin{equation}\label{lf:jt}
j(t):=\min_{(\vy_t,p_t)\in V_g(\oo_t)\times
Q(\oo_t)}\quad\max_{(\vvv_t,q_t)\in V_0(\oo_t)\times
Q(\oo_t)}G(\oo_t,\vy_t,p_t,\vvv_t,q_t),
\end{equation}
$(\vy_t,p_t)$ and $(\vvv_t,q_t)$ satisfy \eqref{state:weak} and
\eqref{adj:weak} on the perturbed domain $\oo_t$, respectively.

Unfortunately, the Sobolev space $V_g(\oo_t)$, $V_0(\oo_t)$, and
$Q(\oo_t)$ depend on the parameter $t$, so we need to introduce the
so-called {function space parametrization} technique which consists
in transporting the different quantities (such as, a cost
functional) defined on the variable domain $\oo_t$ back into the
reference domain $\oo$ which does not depend on the perturbation
parameter $t$. Thus we can use differential calculus since the
functionals involved are defined in a fixed domain $\oo$ with
respect to the parameter $t$.

To do this, we define the following parametrizations
\begin{eqnarray*}
  V_g(\oo_t)&=&\{\vy\circ T_t^{-1}:\;\vy\in V_g(\oo)\};\\
  V_0(\oo_t)&=&\{\vvv\circ T_t^{-1}:\;\vvv\in V_0(\oo)\};\\
  Q(\oo_t)&=&\{p\circ T_t^{-1}:\;p\in Q(\oo)\}.
\end{eqnarray*}
where "$\circ$" denotes the composition of the two maps.

Note that since $T_t$ and $T_t^{-1}$ are diffeomorphisms, these
parametrizations can not change the value of the saddle point. We
can rewrite (\ref{lf:jt}) as
\begin{equation}\label{nsfsp:newsaddlep}
j(t)= \min_{(\vy,p)\in V_g(\oo)\times Q(\oo)}\quad\max_{(\vvv,q)\in
V_0(\oo)\times Q(\oo)}G(\oo_t,\vy\circ T_t^{-1},p\circ
T_t^{-1},\vvv\circ T_t^{-1},q\circ T_t^{-1}).
\end{equation}
where the Lagrangian
$$G(\oo_t,\vy\circ T_t^{-1},p\circ
T_t^{-1},\vvv\circ T_t^{-1},q\circ T_t^{-1})=I_1(t)+I_2(t)+I_3(t)$$
with
$$I_1(t):=2\nu\int_{\oo_t}\seminorm{{\varepsilon}(\vy\circ T_t^{-1})}^2\dx,$$
\begin{multline*}
    I_2(t):=-\int_{\oo_t}[2\nu\ve(\vvv\circ T_t^{-1}):\ve(\vy\circ
  T_t^{-1})+\md(\vy\circ T_t^{-1})\cdot(\vy\circ
  T_t^{-1})\cdot(\vvv\circ T_t^{-1})\\-(p\circ
  T_t^{-1})\mdiv(\vvv\circ T_t^{-1})-(\vy\circ T_t^{-1})\cdot\n
  (q\circ T_t^{-1})]\dx,
\end{multline*}
and
\begin{equation*}
  I_3(t):=\int_{\G_d}\vec h\cdot\vvv\ds.
\end{equation*}
Now we introduce the theorem concerning on the differentiability of
a saddle point (or a minimax). To begin with, some notations are
given as follows.

 Define a functional
$$\mg : [0,\tau]\times X\times Y\rightarrow\mathbb{R}$$
with $\tau>0$, and $X,Y$ are the two topological spaces.

 For any
$t\in [0,\tau]$, define $g(t)=\inf_{x\in X}\sup_{y\in Y}\mg(t,x,y)$
and the sets
\begin{eqnarray*}
  &X(t)=\{x^t\in X:g(t)=\sup_{y\in Y}\mg(t,x^t,y)\}\\
  &Y(t,x)=\{y^t\in Y:\mg(t,x,y^t)=\sup_{y\in Y}\mg(t,x,y)\}
\end{eqnarray*}
Similarly, we can define the dual functional $h(t)=\sup_{y\in
Y}\inf_{x\in X}\mg(t,x,y)$ and the corresponding sets
\begin{eqnarray*}
 & Y(t)=\{y^t\in Y:h(t)=\inf_{x\in X}\mg(t,x,y^t)\}\\
  &X(t,y)=\{x^t\in X:\mg(t,x^t,y)=\inf_{x\in X}\mg(t,x,y)\}
\end{eqnarray*}
Furthermore, we introduce the set of saddle points
$$S(t)=\{(x,y)\in X\times Y: g(t)=\mg(t,x,y)=h(t)\}$$
Now we can introduce the following theorem (see \cite{correa} or
page 427 of \cite{delfour}):
\begin{theorem}\label{fsp:correa}
 Assume that the following hypothesis hold:
 \begin{itemize}
    \item [(H1)]$S(t)\neq\emptyset,\;t\in [0,\tau];$
    \item [(H2)]The partial derivative $\p_t\mg(t,x,y)$ exists in
    $[0,\tau]$ for all $$(x,y)\in \left[\underset{{t\in [0,\tau]}}{\bigcup}X(t)\times Y(0)\right]\bigcup\left[X(0)\times\underset{{t\in [0,\tau]}}{\bigcup}Y(t)\right];$$
    \item [(H3)]There exists a topology $\mt_X$ on $X$ such that for
    any sequence $\{t_n:t_n\in [0,\tau]\}$ with
    $\lim\limits_{n\nearrow\infty}t_n=0$, there exists $x^0\in X(0)$ and a subsequence
    $\{t_{n_k}\}$, and for each $k\geq 1,$ there exists $x_{n_k}\in
    X(t_{n_k})$ such that
    \begin{enumerate}
        \item [(i)]$\lim\limits_{n\nearrow\infty}x_{n_k}=x^0$ in the
        $\mt_X$ topology,
        \item [(ii)]$\liminf\limits_{t\searrow 0\atop k\nearrow\infty}\p_t\mg(t,x_{n_k},y)\geq\p_t\mg(0,x^0,y),
        \quad \forall y\in Y(0);$
    \end{enumerate}
    \item [(H4)]There exists a topology $\mt_Y$ on $Y$ such that for
    any sequence $\{t_n:t_n\in [0,\tau]\}$ with
    $\lim\limits_{n\nearrow\infty}t_n=0$, there exists $y^0\in Y(0)$ and a subsequence
    $\{t_{n_k}\}$, and for each $k\geq 1,$ there exists $y_{n_k}\in
    Y(t_{n_k})$ such that
    \begin{enumerate}
        \item [(i)]$\lim\limits_{n\nearrow\infty}y_{n_k}=y^0$ in the
        $\mt_Y$ topology,
        \item [(ii)]$\limsup\limits_{t\searrow 0\atop k\nearrow\infty}
        \p_t\mg(t,x,y_{n_k})\leq\p_t\mg(0,x,y^0),\quad \forall x\in X(0).$
    \end{enumerate}
 \end{itemize}
 Then there exists $(x^0,y^0)\in X(0)\times Y(0)$ such that
 \begin{multline}
   \diff g(0)=\lim_{t\searrow 0}\frac{g(t)-g(0)}{t}\\=\inf_{x\in
   X(0)}\sup_{y\in Y(0)}\p_t \mg(0,x,y)=\p_t\mg(0,x^0,y^0)=\sup_{y\in Y(0)}\inf_{x\in
   X(0)}\p_t \mg(0,x,y)
 \end{multline}
 This means that $(x^0,y^0)\in X(0)\times Y(0)$ is a saddle point of
 $\p_t\mg(0,x,y)$.
\end{theorem}
Following Theorem \ref{fsp:correa}, we need to differentiate the
perturbed Lagrange functional $G(\oo_t,\vy\circ T_t^{-1},p\circ
T_t^{-1},\vvv\circ T_t^{-1},q\circ T_t^{-1})$.

To perform the differentiation, we introduce the following Hadamard
formula\cite{ha07}
\begin{equation}\label{hadamard}
 \frac{\diff{}}{\diff t}\int_{\oo_t}F(t,x)\dx=\int_{\oo_t}
 \frac{\p F}{\p t}(t,x)\dx+\int_{\p\oo_t} F(t,x)\,\vec
 V\cdot\vn_t\ds,
 \end{equation}
 for a sufficiently smooth functional
$F:[0,\tau]\times\rn\rightarrow\mathbb{R}$.

By Hadamard formula (\ref{hadamard}), we get
$$\p_t G(\oo_t,\vy\circ T_t^{-1},p\circ
T_t^{-1},\vvv\circ T_t^{-1},q\circ
T_t^{-1})=I'_1(0)+I'_2(0)+I'_3(0)+I'_4(0),$$ where
\begin{equation}\label{i1}
I'_1(0)=4\nu\int_\oo
\ve(\vy):\ve(-\md\vy\cdot\vv)\dx+2\nu\int_{\G_s}\seminorm{\ve(\vy)}^2\vv_n\ds;
\end{equation}
\begin{multline}\label{i2}
  I'_2(0)=-\int_\oo[2\nu\ve(-\md\vy\cdot\vv)\cdot\ve(\vvv)+2\nu\ve(\vy)\cdot\ve(-\md\vvv\cdot\vv)
  +\md\vy\cdot\vy\cdot(-\md\vy\cdot\vv)\\+\md(-\md\vy\cdot\vv)\cdot\vy\cdot\vvv+\md\vy\cdot(-\md\vy\cdot\vv)\cdot\vvv
   -p\mdiv(-\md\vvv\cdot\vv)\\-\mdiv(-\md\vy\cdot\vv)q
  -\mdiv\vy(-\n q\cdot\vv)-(-\n p\cdot\vv)\mdiv\vvv]\dx\\
  +\int_{\G_s}(-2\nu\ve(\vy):\ve(\vvv)-\md\vy\cdot\vy\cdot\vvv+p\mdiv\vvv+\mdiv\vy
  q)\vv_n\ds;
\end{multline}
and
$
  I'_3(0)=0.
$

To simplify \eqref{i1} and \eqref{i2}, we introduce the following
lemma.
\begin{lemma}\label{lem:a}
  If two vector functions $\vy$ and $\vvv$ vanish on the boundary
  $\G_s$ and $\mdiv\vy=\mdiv\vvv=0$ in $\oo$, the
  following identities
  \begin{eqnarray}
   \label{lem:a1}  &  \md\vy\cdot\vv\cdot\vn=(\md\vy\cdot\vn\cdot\vn)\vv_n=\mdiv\vy\vv_n;\\
  \label{lem:a2}   &\ve(\vy):\ve(\vvv)=\ve(\vy):(\ve(\vvv)\cdot(\vn\otimes\vn))=(\ve(\vy)\cdot\vn)\cdot(\ve(\vvv)\cdot\vn);\\
  \label{lem:a3}  &(\ve(\vy)\cdot\vn)\cdot(\md\vvv\cdot\vv)=(\ve(\vy)\cdot\vn)\cdot(\md\vvv\cdot\vn)\vv_n
  =(\ve(\vy)\cdot\vn)\cdot(\ve(\vvv)\cdot\vn)\vv_n
  \end{eqnarray}
  hold on the boundary $\G_s$, where the tensor product $\vn\otimes\vn:=\sum_{i,j=1}^Nn_i n_j.$
\end{lemma}
Using Lemma \ref{lem:green}, for \eqref{i1} we have
\begin{equation*}
  \label{i1:a0}
  I'_1(0)=-2\nu\int_\oo\Delta\vy\cdot(-\md\vy\cdot\vv)\dx+4\nu\int_{\G_s}
  (\ve(\vy)\cdot\vn)\cdot(-\md\vy\cdot\vv)\ds+2\nu\int_{\G_s}\seminorm{\ve(\vy)}^2\vv_n\ds.
\end{equation*}
By the identities \eqref{lem:a2} and \eqref{lem:a3}, we further get
\begin{equation}\label{i1:a}
  I'_1(0)=-2\nu\int_\oo\Delta\vy\cdot(-\md\vy\cdot\vv)\dx-2\nu\int_{\G_s}\seminorm{\ve(\vy)}^2\vv_n\ds.
\end{equation}
Employing Lemma \ref{lem:green} and
$\vy|_{\G_s}=\vv|_{\G_w\cup\G_u\cup\G_d}=0$, \eqref{i2} can be
rewritten as
\begin{multline}\label{i2:a}
  I'_2(0)=\int_\oo [(\nu\Delta\vy-\md\vy\cdot\vy-\n
  p)\cdot(-\md\vvv\cdot\vv)+\mdiv\vy(-\n q\cdot\vv)]\dx\\
  +\int_\oo[(\nu\Delta\vvv+\md\vvv\cdot\vy-{}^*\md\vy\cdot\vvv-\n
  q)\cdot(-\md\vy\cdot\vv)+\mdiv\vvv(-\n p\cdot\vv)]\dx\\
  -\int_{\G_s}[\sigma(\vy,p)\cdot\vn\cdot(-\md\vvv\cdot\vv)
  +\sigma(\vvv,q)\cdot\vn\cdot(-\md\vy\cdot\vv)]\ds\\
  -\int_{\G_s}[2\nu\ve(\vy):\ve(\vvv)+\md\vy\cdot\vy\cdot\vvv-p\mdiv\vvv-\mdiv\vy
  q]\vv_n\ds.
\end{multline}
Since $(\vy,p)$ and $(\vvv,q)$ satisfy \eqref{nsp:a}\eqref{nsp:b}
and \eqref{adj:b}\eqref{adj:a} respectively, \eqref{i2:a} reduces to
\begin{multline}\label{i2:b}
   I'_2(0)=2\nu\int_\oo\Delta\vy\cdot(-\md\vy\cdot\vv)\dx\\-\int_{\G_s}[\sigma(\vy,p)\cdot\vn\cdot(-\md\vvv\cdot\vv)
  +\sigma(\vvv,q)\cdot\vn\cdot(-\md\vy\cdot\vv)+2\nu\ve(\vy):\ve(\vvv)\vv_n]\ds.
\end{multline}
On the boundary $\G_s$, we can deduce that
\begin{equation*}
\begin{array}{@{\hspace*{2cm}}ll@{\hspace{2cm}}r}
 \lefteqn{ -\sigma(\vy,p)\cdot\vn\cdot(-\md\vvv\cdot\vv)
  -\sigma(\vvv,q)\cdot\vn\cdot(-\md\vy\cdot\vv)\hspace*{3cm}}\\
  =&2\nu[\ve(\vy)\cdot\vn\cdot(\md\vvv\cdot\vv)+\ve(\vvv)\cdot\vn\cdot(\md\vy\cdot\vv)]&(\textrm{by}\;\eqref{lem:a1})\\
  =&4\nu(\ve(\vy)\cdot\vn)\cdot(\ve(\vvv)\cdot\vn)\vv_n &(\textrm{by}\;\eqref{lem:a3})\\
  =&4\nu\ve(\vy):\ve(\vvv)\vv_n.
  &(\textrm{by}\;\eqref{lem:a2})
\end{array}
\end{equation*}
Therefore, \eqref{i2:b} becomes
\begin{equation}\label{i2:c}
  I'_2(0)=2\nu\int_\oo\Delta\vy\cdot(-\md\vy\cdot\vv)\dx+2\nu\int_{\G_s}\ve(\vy):\ve(\vvv)\vv_n\ds.
\end{equation}
Adding \eqref{i1:a} and \eqref{i2:c} together, we finally obtain the
boundary expression for the Eulerian derivative of $J(\oo)$,
\begin{equation}
  \diff J(\oo;\vec V)=2\nu\int_{\G_s}\left[\ve(\vy):\ve(\vvv)-\seminorm{\ve(\vy)}^2\right]\vv_n\ds,
\end{equation}
Since the mapping $\vec V\mapsto \diff J(\oo;\vec V)$ is linear and
continuous, we get the expression for the shape gradient
\begin{equation}
   \n J=2\nu[\ve(\vy):\ve(\vvv)-\seminorm{\ve(\vy)}^2]\vn
\end{equation} by (\ref{pri:shaped}).
\section{Finite element approximations and numerical Simulation}

\subsection{Discretization of the optimization problem}
We suppose that $\oo$ is a bounded polygonal domain of
$\mathbb{R}^2$ and only consider the conforming finite element
approximations. Let $X_h\subset H^1(\oo)^N$ and $S_h\subset
L^2(\oo)$ be two families of finite dimensional subspaces
parameterized by $h$ which tends to zero. We also define
\begin{eqnarray*}
  V_{gh}&:=&\{\vu_h\in X_h: \;\vu_h=0\mbox{ on }\G_w\cup\G_s,\;\vu_h=\vg\mbox{ on }\G_u\},\\
V_{0h}&:=&\{\vu_h\in X_h: \;\vu_h=0\mbox{ on }\G_w\cup
\G_u\cup \G_s\},\\
 Q_h&:=&\left\{p_h\in S_h:\;\int_\oo p_h\dx=0 \;(\mbox{ if
meas}(\G_d)=0)\right\}.
\end{eqnarray*}
Besides, the following assumptions are supposed to hold.
\begin{itemize}
  \item [(HA1)]There exists $C>0$ such that for $0\leq m\leq l$,
  $$\inf_{\vvv_h\in V_{gh}}\norm{\vvv_h-\vvv}_1\leq C
  h^m\norm{\vvv}_{m+1},\qquad \forall \vvv\in H^{m+1}(\oo)^N\cap
  V_g(\oo);$$
  \item [(HA2)]There exists $C>0$ such that for $0\leq m\leq l'$,
  $$\inf_{q_h\in Q_h}\norm{q_h-q}_0\leq Ch^m\norm{q}_m,\quad
  \forall q\in H^m(\oo)\cap Q(\oo);$$
  \item [(HA3)]The Ladyzhenskaya-Brezzi-Babuska inf-sup condition is
  verified, i.e., there exists $C>0$, such that
  $$ \inf_{0\neq q_h\in Q_h}\sup_{0\neq \vvv_h\in V_h}\frac{\int_\oo q_h\mdiv\vvv_h\dx}{\norm{\vvv_h}_1\norm{q_h}_0}\geq
C,\qquad V_h=V_{gh}\mbox{ or }V_{0h}.$$
\end{itemize}
The Galerkin finite element approximations of the state system
(\ref{state:weak}) and adjoint state system (\ref{adj:weak}) in
mixed form are as follows
\begin{equation}\label{state:weak2}
  \left\{
  \begin{array}
    {ll}
    &\mbox{seek } (\vy_h,p_h)\in V_{gh} \times Q_h\mbox{ such
  that}\;\forall (\vvv_h,q_h)\in V_{0h}\times Q_h,\\[4pt]
  &   \int_{\oo}[2\nu\varepsilon(\vy_h):\varepsilon(\vvv_h)
     +\md\vy_h\cdot\vy_h\cdot\vvv_h-p_h\mdiv\vvv_h]\dx=\int_{\G_d}\vec
     h\cdot\vvv_h\ds,\\[4pt]
 &\int_{\oo}\mdiv\vy_h q_h\dx=0,
  \end{array}
  \right.
\end{equation}
and
\begin{equation}\label{adj:weak2}
  \left\{
  \begin{array}
    {ll}
    &\mbox{seek } (\vvv_h,q_h)\in V_{0h}\times Q_h\mbox{ such
  that}\;\forall (\vphi_h,\pi_h)\in V_{0h}\times Q_h,\\[5pt]
  &   \int_{\oo}[2\nu\varepsilon(\vvv_h):\varepsilon(\vphi_h)
     +\md\vphi_h\cdot\vy_h\cdot\vvv_h+\md\vy_h\cdot\vphi_h\cdot\vvv_h-q_h\mdiv\vphi_h]\dx\\
    &\hspace{5cm} =4\nu\int_\oo\ve(\vy_h):\ve(\vphi_h)\dx,\\[4pt]
 &\int_{\oo}\mdiv\vvv_h \pi_h\dx=0.
  \end{array}
  \right.
\end{equation}
We also have the discrete cost functional
\begin{equation}
  J_{h}(\oo)=2\int_\oo\seminorm{\ve(\vy_h)}^2\dx,
\end{equation}
and the discrete shape gradient
\begin{equation}
  \n J_h=2\nu[\ve(\vy_h):\ve(\vvv_h)-\seminorm{\ve(\vy_h)}^2]\vn
\end{equation}
Finally for completeness, we state the following theorem (see
\cite{girault86}).
\begin{theorem}
  Assume that the hypotheses {(HA1)}, (HA2) and (HA3) hold. Let
  $$\{(\lambda,(\vy(\lambda),\lambda p(\lambda))); \lambda=1/\nu\in
  \Lambda,\;\Lambda \mbox{ is a connected subsect of } \mathbb{R}^+\}$$ be a branch of nonsingular solutions of the
  state system \eqref{state:weak}. Then there exists a neighborhood
  $\mathcal{O}$ of the origin in $V_g(\oo)\times Q(\oo)$ and for
  $h\leq h_0$ sufficiently small a unique $C^\infty$ branch $\{(\lambda,(\vy_h(\lambda),\lambda p_h(\lambda))); \lambda\in
  \Lambda\}$ of nonsingular solutions of problem \eqref{state:weak2}
  such that
 \begin{equation*}
   \lim_{h\rightarrow
   0}\sup_{\lambda\in\Lambda}\{\norm{\vy_h(\lambda)-\vy(\lambda)}_2+\norm{p_h(\lambda)-p(\lambda)}_1\}=0.
  \end{equation*}
  In addition, for the adjoint state system \eqref{adj:weak} and its
  discrete form \eqref{adj:weak2}, we have the similar convergence
  result.
\end{theorem}

\subsection{A gradient type algorithm}
For the minimization problem \eqref{nsdrag:cost}, we rather work
with the unconstrained minimization problem
\begin{equation}\label{nsdrag:cost2}
  \min_{\oo\in \mathbb{R}^2}G(\oo)=J(\oo)+lV(\oo),
\end{equation}
where $V(\oo):=\int_\oo\dx$ and $l$ is a positive Lagrange
multiplier. The Eulerian derivative of $G(\oo)$ is
\begin{equation*}
  \diff G(\oo;\vec V)=\int_{\G_s}\n G\cdot\vv \ds,
\end{equation*}
where the shape gradient $\n
G:=[2\nu\ve(\vy):\ve(\vvv)-2\nu\seminorm{\ve(\vy)}^2+l]\vn$.
Ignoring regularization, a descent direction is found by defining $
  \vec V=-h_k\n G$, and then we can update the shape $\oo$ as $
\oo_k=(\mathrm{I}+h_k\vec V)\oo $, where $h_k$ is a descent step at
$k$-th iteration.

However, in this article in order to avoid boundary oscillations
(and irregular shapes) and due to the fact that the gradient type
algorithm produces shape variations which have less regularity than
the original parametrization, we change the scalar product with
respect to which we compute a descent direction, for instance,
$H^1(\oo)^2$. In this case, the descent direction is the unique
element $\vec d\in H^1(\oo)^2$ of the problem
\begin{equation}\label{reg}
  \left\{
  \begin{array}{ll}
  -\Delta \vec d+\vec d=0 \quad&\mbox{in }\oo,\\
  \vec d=0,&\mbox{on }\G_u\cup\G_d\cup\G_w,\\
  \md\vec d\cdot\vn=-\n G\;&\mbox{on }\G_s.
  \end{array}
\right.
\end{equation}
 To better understand the necessity of projection or smoother due
  to the loss of regularity, we give the following remark.
\begin{rem}We give a simple example to illustrate the loss of
regularity. We suppose that the cost functional is a quadratic
functional: $J(x)=(Ax-b)^2$ with $x\in H^1(\oo)$, $A\in H^{-1}(\oo)$
and $b\in L^2(\oo)$. The gradient $\n J=2(Ax-b)A\in H^{-1}(\oo)$ has
less regularity than $x$. Then any variation using $\n J$ as the
descent direction will have less regularity than $x$, therefore we
need to project into $H^1(\oo)$. We refer the readers to see
B.Mohammadi $\&$ O.Pironneau \cite{piron01} and G.Dogan
{et.al.}\cite{dogan06} for further discussion on regularity.
\end{rem}
The resulting algorithm can be summarized as follows:
\begin{itemize}
    \item [(1)] Choose an initial shape $\oo_0$, an initial step $h_0$ and a Lagrange multiplier $l_0$;
      \item [(2)] Compute the state system \eqref{state:weak} and the adjoint system \eqref{adj:weak}, then
    we can evaluate the descent direction $\vec d_k$ by using (\ref{reg})
    with $\oo=\oo_k$ and $l=l_k$;
    \item[(3)] Set $\oo_{k+1}=(\mathrm{I}+h_k\vec d_k) \,\oo_k$
    and $l_{k+1}=(l_k+l)/2+\epsilon\seminorm{V(\oo_k)-V(\oo)}/V(\oo)$ with a small positive constant $\epsilon$,
    where $l=-\int_{\G_s}\n
    J\ds/\int_{\G_s}\ds$ and $V(\oo)$ is the given area of $\oo$.
\end{itemize}
The choice of the descent step size $h_k$ is not an easy task. Too
big, the algorithm is unstable; too small, the rate of convergence
is insignificant. In order to refresh $h_k$, we compare $h_k$ with
$h_{k-1}$. If $(\vec d_k,\vec d_{k-1})_{H^1}$ is negative, we should
reduce the step; on the other hand, if $\vec d_k$ and $\vec d_{k-1}$
are very close, we increase the step. In addition, if reversed
triangles are appeared when moving the mesh, we also need to reduce
the step.

In our algorithm, we do not choose any stopping criterion. A
classical stopping criterion is to find that whether the shape
gradients in some suitable norm is small enough. However, since we
use the continuous shape gradients, it's hopeless for us to expect
very small gradient norm because of numerical discretization errors.
Instead, we fix the number of iterations. If it is too small, we can
restart it with the previous final shape as the initial shape.

\subsection{Numerical results}
In all computations, the finite element discretization is effected
using the $P_{1}$bubble--$P_1$ pair of finite element spaces on a
triangular mesh, i.e., we choose the following velocity space $X_h$
and pressure space $S_h$:
\begin{eqnarray*}
  X_h&=&\{\vy_h\in (C^0(\bar{\oo}))^2: \vy_h|_T\in (P^*_{1T})^2, \forall
  T\in\mathcal{T}_h\}\\
  S_h&=&\{p_h\in C^0(\bar{\oo}): p_h|_T\in P_1, \forall T\in
  \mathcal{T}_h\},
\end{eqnarray*}
where $\mathcal{T}_h$ denotes a standard finite element
triangulation of $\oo$, $P_k$ the space of the polynomials in two
variables of degree $\leq k$ and $P^*_{1T}$ the subspace of $P_3$
defined by
\begin{multline*}
  P^*_{1T}=\{q: q=q_1+\lambda\phi_T, \mbox{ with } q_1\in P_1, \lambda\in\mathbb{R}\;\mbox{and}\\
                 \phi_T\in P_3, \phi_T=0\mbox{ on }\p T,\; \phi_T(G_T)=1\mbox{ with }G_T\mbox{ is the centroid of
                 }T\}.
\end{multline*}
Notice that a function like $\phi_T$ is usually called a bubble
function.

The mesh is performed by a Delaunay-Voronoi mesh generator (see
\cite{piron01}) and during the shape deformation, we utilize the a
metric-based anisotropic mesh adaptation technique where the metric
can be computed automatically from the Hessian of a solution. We run
the programs on a home PC with Intel Pentium 4 CPU 2.8 GHz and 1GB
memory.
\subsubsection{Test case 1: cannula shape optimization in Stokes flow}
We consider the shape optimization of a two-dimensional inflow
cannula of a circulatory assist device in the biomedical
applications. The geometry of the cannula $\oo$ is depicted in the
left picture of \autoref{fig:cannula}. The boundary conditions for
the problem are traction-free at the exit $\G_d$, no-slip at all
curved walls $\G_s$, and a specified parabolic inlet velocity
$\vg(0,y)=((y-2)(2.35-y),0)^T$.
\begin{figure}[!htbp]
\renewcommand{\captionlabelfont}{\small}
\setcaptionwidth{5.5in}
\begin{minipage}[b]{0.45\textwidth}
  \centering
 {\includegraphics[width=2.45in]{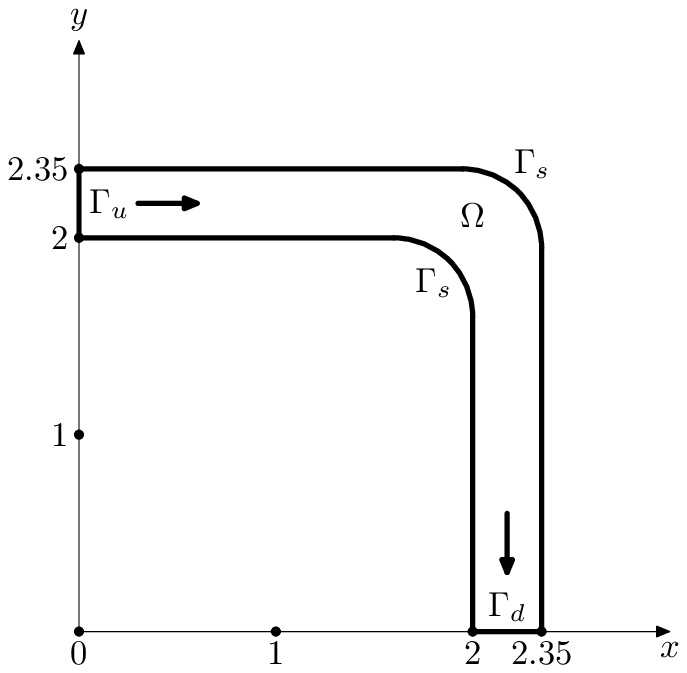}}\\
  \end{minipage}
\begin{minipage}[b]{0.45\textwidth}
  \centering
{\includegraphics[width=1.95in]{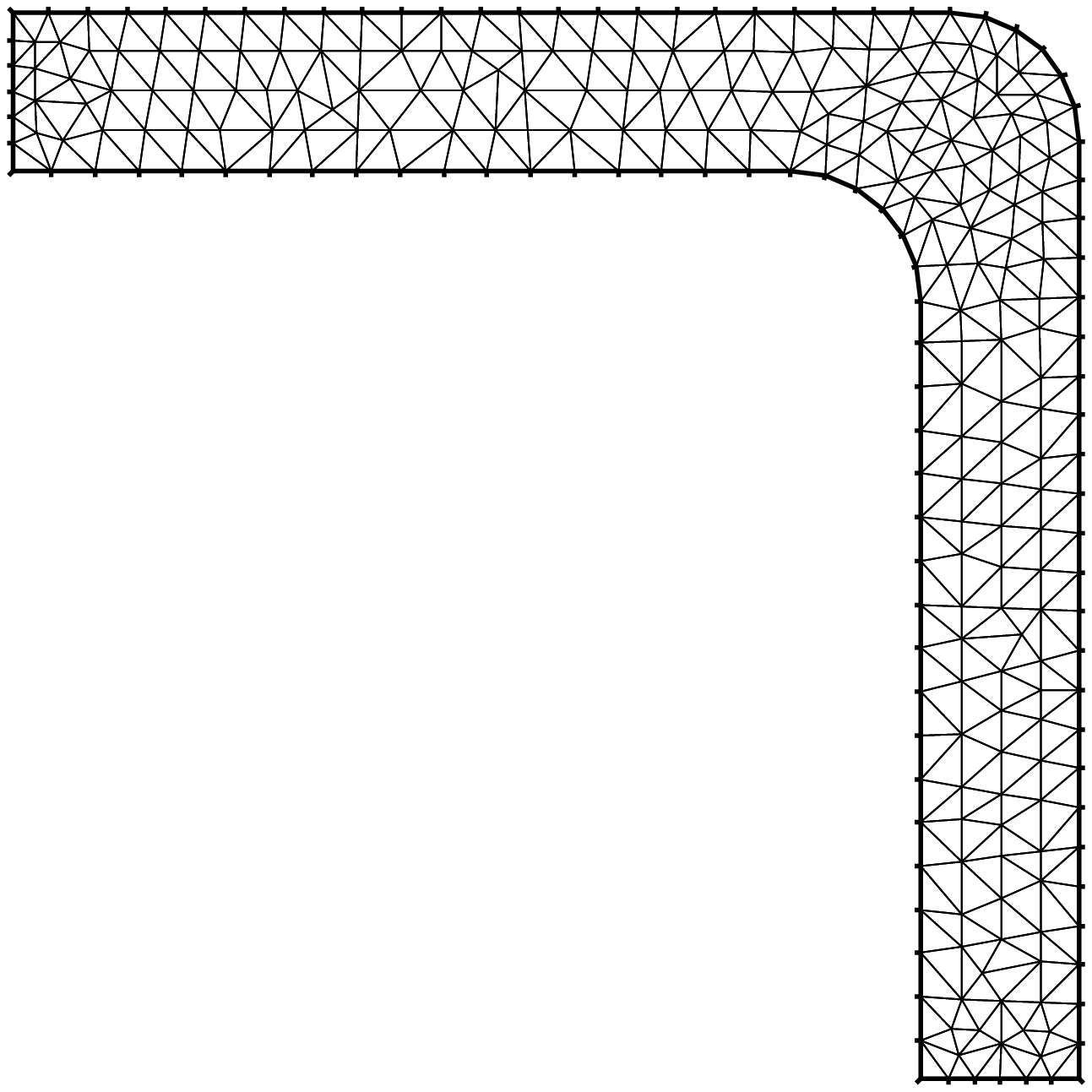}}\\
  \end{minipage}
  \caption{The analytic domain and the finite element mesh of the cannula.\label{fig:cannula}}
  \end{figure}

In this test case, we present results for two different Reynolds
numbers 0.1 and 0.01, defined as
$\mathrm{Re}=d\seminorm{\vy_{\mathrm{m}}}/\nu$, where
$\vy_{\mathrm{m}}$ is the maximum velocity at the inlet $\G_u$ and
$d=0.35$ is the diameter of the cannula. The domain is discretized
using 448 triangular elements and 279 nodes. Since the inertial term
in \eqref{nsp:a} can be neglected when $\mathrm{Re}=0.1$ or $0.01$,
we can say that the blood flow in the cannula was governed by the
Stokes equations approximately.
  \begin{figure}[h]
\renewcommand{\captionlabelfont}{\small}
\setcaptionwidth{5.5in}
\begin{minipage}[b]{0.45\textwidth}
  \centering
 {\includegraphics[width=2.5in]{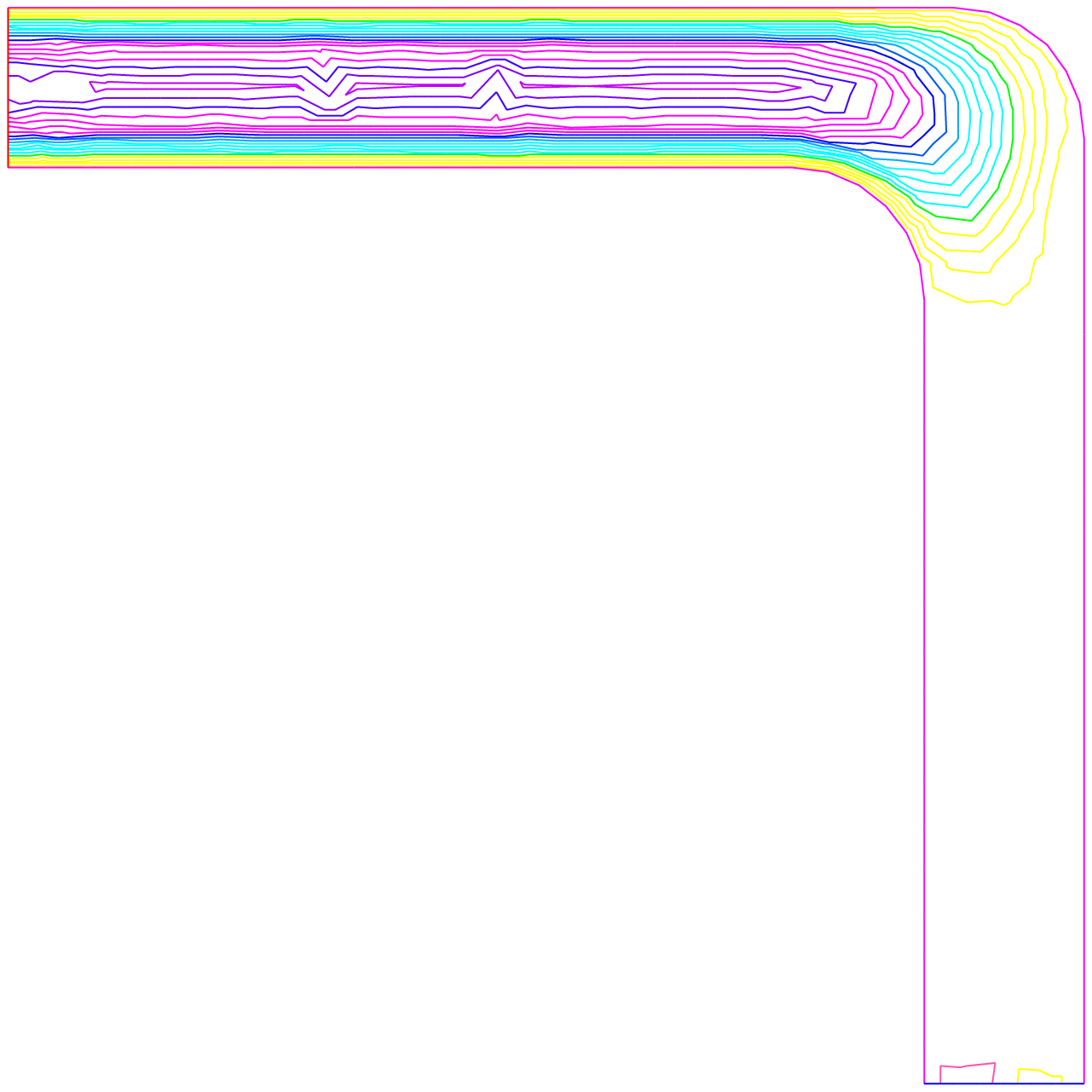}}\\
  \end{minipage}
\begin{minipage}[b]{0.45\textwidth}
  \centering
{\includegraphics[width=2.5in]{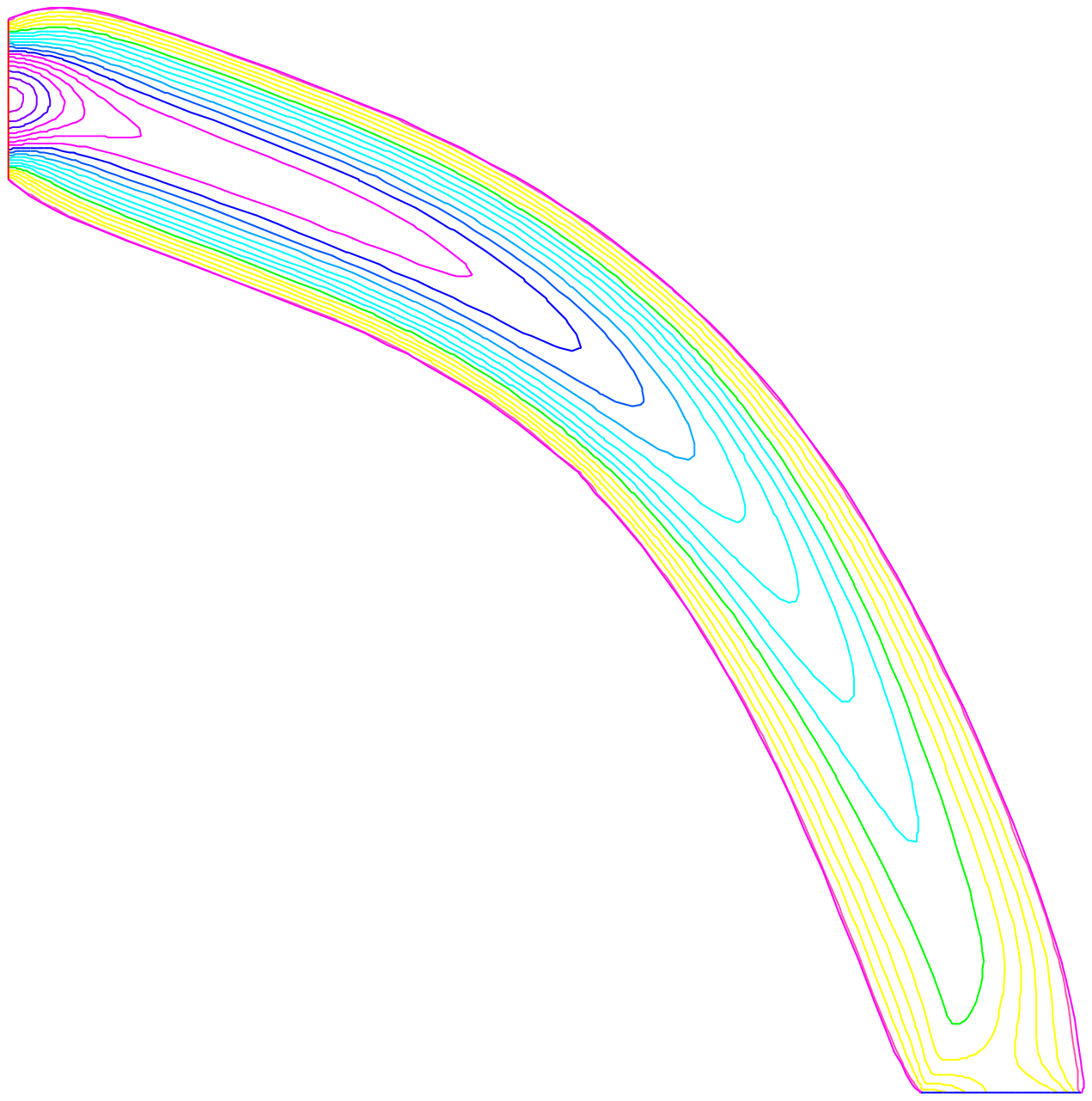}}\\
  \end{minipage}
  \caption{The initial and optimal cannula shapes with
  $\mathrm{Re}=0.1$.\label{fig:can1}}
  \end{figure}
\begin{figure}[!htbp]
\renewcommand{\captionlabelfont}{\small}
\setcaptionwidth{5.5in}
\begin{minipage}[b]{0.45\textwidth}
  \centering
 {\includegraphics[width=2.5in]{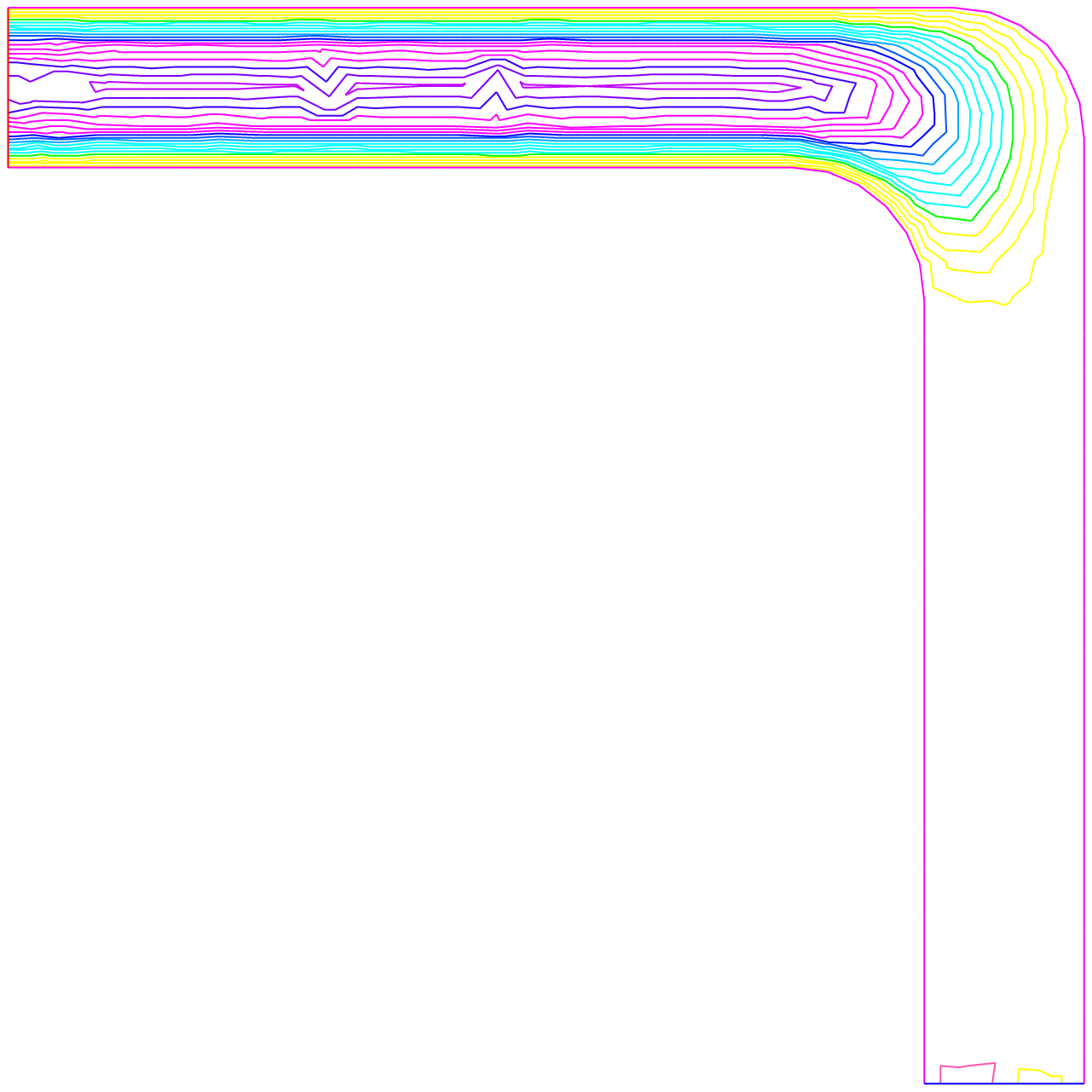}}\\
  \end{minipage}
\begin{minipage}[b]{0.45\textwidth}
  \centering
{\includegraphics[width=2.5in]{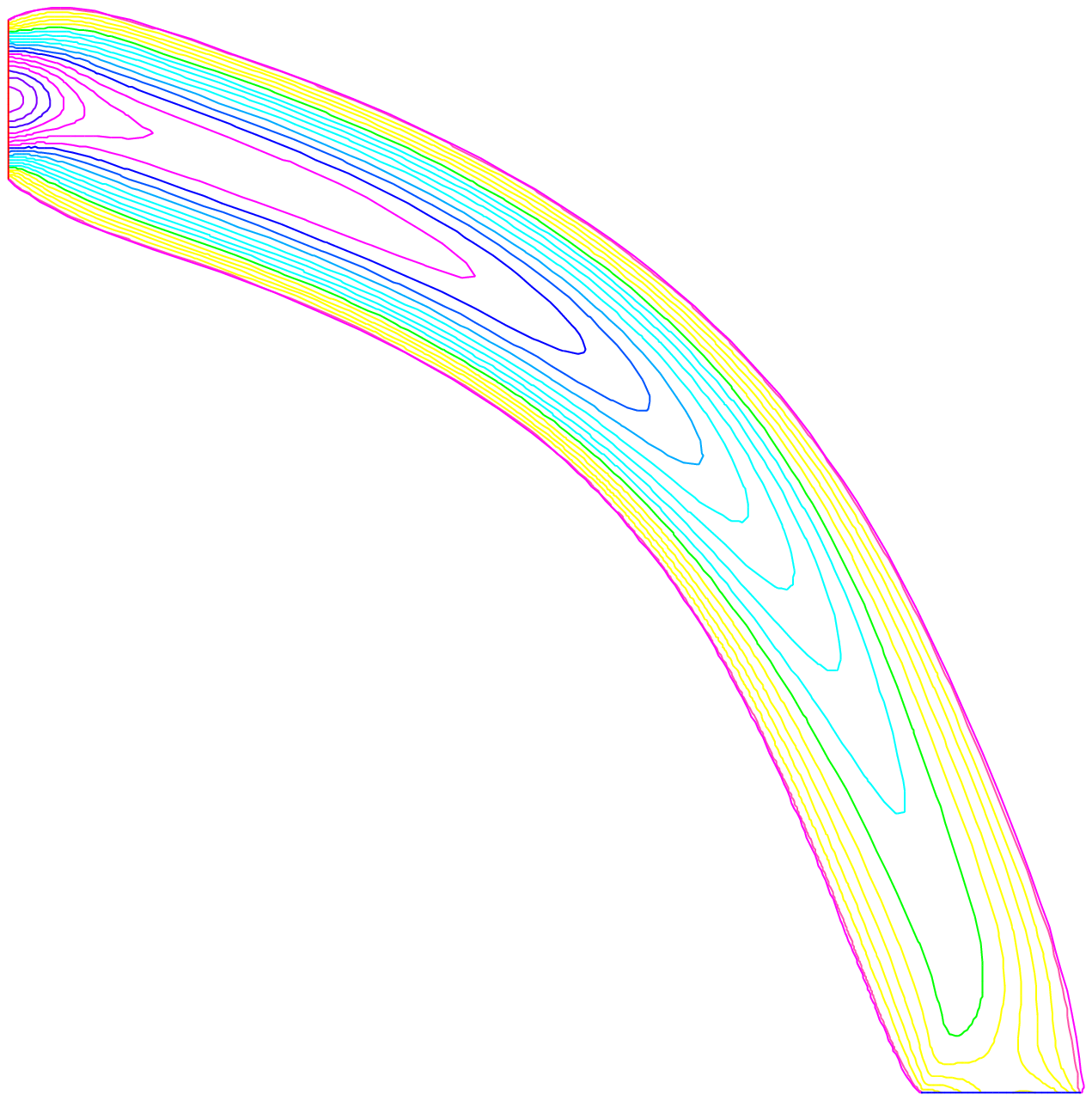}}\\
  \end{minipage}
  \caption{The initial and optimal cannula shapes with
  $\mathrm{Re}=0.01$.\label{fig:can2}}
  \end{figure}

The distributions of the horizontal velocity for the initial and
optimal shapes with $\mathrm{Re}=0.1,0.01$ are shown in
\autoref{fig:can1} and \autoref{fig:can2}. It is clear that shape
optimization has removed the sharp bend in the initial configuration
of the cannula.

The optimization process gave a 55.4975\% reduction in the
dissipated energy with $\mathrm{Re}=0.1$, and a 55.0392\% reduction
in the dissipated energy with $\mathrm{Re}=0.01$.

\subsubsection{Test case 2: Optimization of a solid body in the
Navier--Stokes flow}
 As a second test case, we consider the isolated
body problem. The schematic geometry of the fluid domain is
described in \autoref{fig:domain}, corresponding to an external flow
around a solid body $S$. We reduce the problem to a bounded domain
$D$ by introducing an artificial boundary $\p
D:=\G_u\cup\G_d\cup\G_w$ which has to be taken sufficiently far from
$S$ so that the corresponding flow is a good approximation of the
unbounded external flow around $S$ and $\oo:=D\backslash\bar{S}$ is
the effective domain. In addition, the boundary $\G_s:=\p S$ is to
be optimized.
\begin{figure}[!htbp]
\renewcommand{\captionlabelfont}{\small}
\setcaptionwidth{5.5in}\centering
 \includegraphics[width=2.5in]{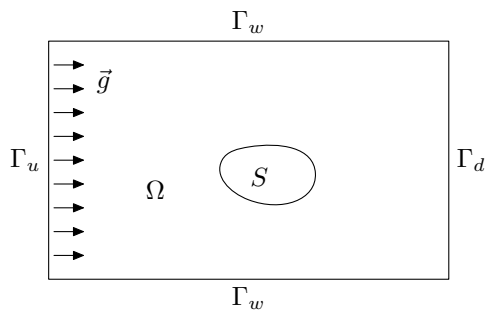}
 \caption{External flow around a solid body $S$.\label{fig:domain}}
  \end{figure}

We choose $D$ to be a rectangle $(-0.5,1.5)\times (-0.5,1.5)$
 and $S$ is to be determined in our simulations. The inflow velocity is assumed to be parabolic with a profile
$\vg(-0.5,y)=(0.2y^2-0.05,0)^T,$ while at the outflow boundary
$\G_d$, we impose a traction-free boundary condition ($\vec h=0$).
No-slip boundary condition are imposed at all the other boundaries.
We further define the admissible set
$$\mathcal{O}:=\left\{\oo\subset\mathbb{R}^2:\; \p D \mbox{ is
fixed},\;\mbox{the area }V(\oo)=1.9\right\},$$ which means that the
target volume of $S$ to be optimized is $0.1$.

We choose the initial shape of the body $S$ to be a circle of center
$(0,0)$ with radius $r=0.2$.
  We present results for two different
Reynolds numbers $\mathrm{Re}=40, 200$ defined by
$\mathrm{Re}=2r\seminorm{\vy_{\mathrm{m}}}/\nu$, where
$\vy_{\mathrm{m}}$ is the maximum velocity at the inflow $\G_u$. The
finite element mesh used for the calculations at $\mathrm{Re}=200$
has been shown in \autoref{fig:fem}.

\begin{figure}[!htbp]
\renewcommand{\captionlabelfont}{\small}
\setcaptionwidth{5.5in}
  \centering
{\includegraphics[width=2.6in]{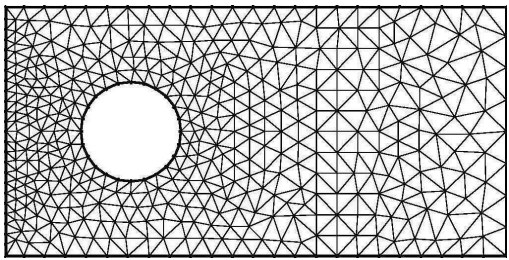}}
  \caption{the finite element mesh.\label{fig:fem}}
  \end{figure}

\autoref{fig0:e} and \autoref{fig0:c} represent the distribution of
the velocity $\vy=(y_1,y_2)^T$ and the pressure $p$ for the initial
shape and the optimal shape in the neighborhood of $S$ for
$\mathrm{Re}=40, 200$ respectively.

We run many iterations in order to show the good convergence and
stability properties of our algorithm, however it is clear that it
has converged in a smaller number of iterations (see
\autoref{fig:cost2} and \autoref{fig:cost3}). For the Reynolds
numbers $40$ and $200$, the total dissipated energy reduced about
$34.47\%$ and $44.46\%$ respectively.
\section{Conclusion}
The minimization problem of total dissipated energy in the two
dimensional Navier--Stokes flow with mixed boundary conditions
involving pressure has been presented. We derived the structure of
shape gradient for the cost functional by function space
parametrization technique without the usual study of the derivative
of the state. Though for the time being this technique lacks from a
rigorous mathematical framework for the Navier--Stokes equations, a
gradient type algorithm is effectively used for the minimization
problem for various Reynolds numbers. Further research is necessary
on efficient implementations for time--dependent Navier--Stokes flow
and much more real problems in the industry.
\begin{figure}[!htbp]
\renewcommand{\captionlabelfont}{\small}
\setcaptionwidth{5.5in}
\begin{minipage}[b]{0.32\textwidth}
  \centering
 {\includegraphics[width=1.6in]{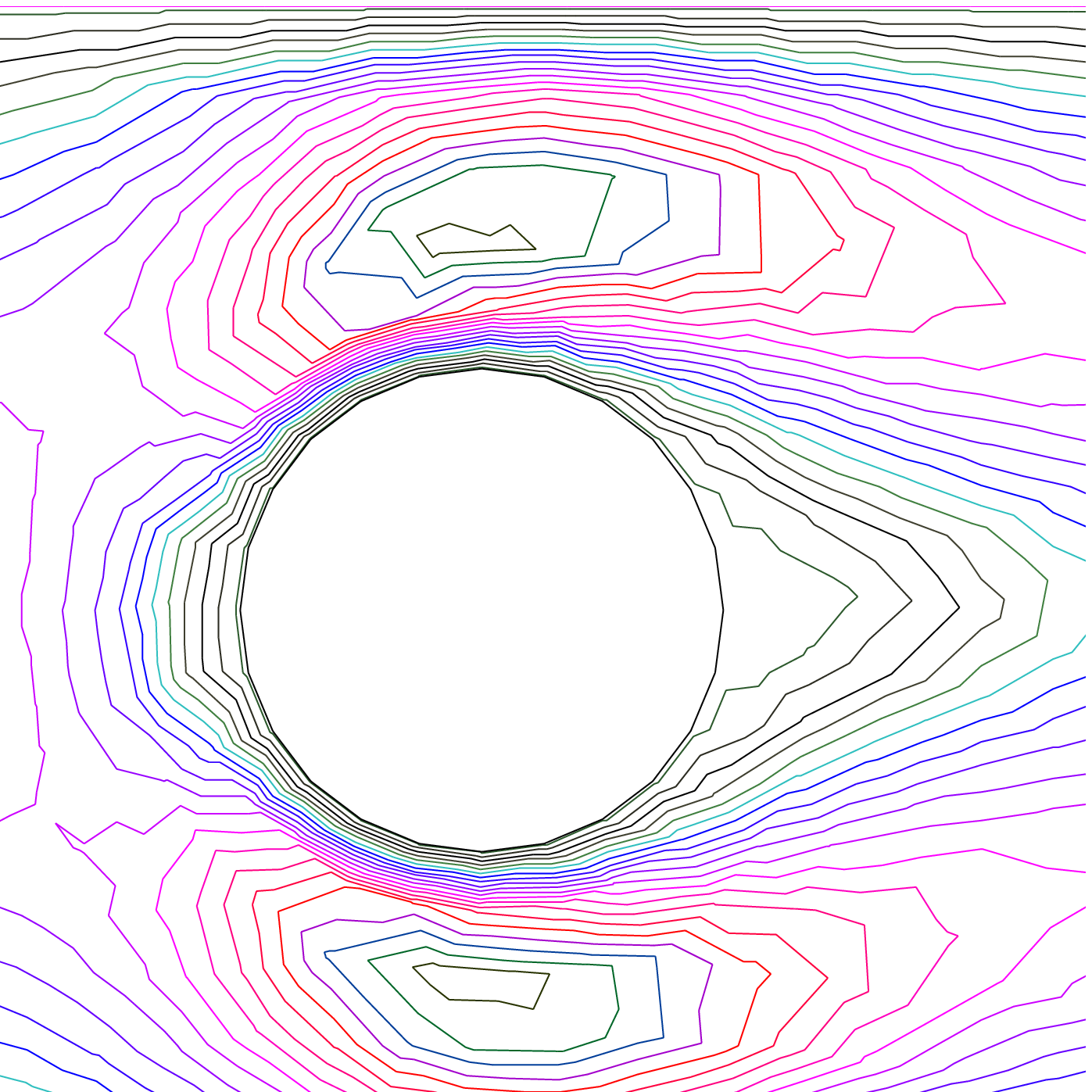}}\\
  {(a) $y_1$ for initial shape.}
  \end{minipage}
\begin{minipage}[b]{0.32\textwidth}
  \centering
{\includegraphics[width=1.6in]{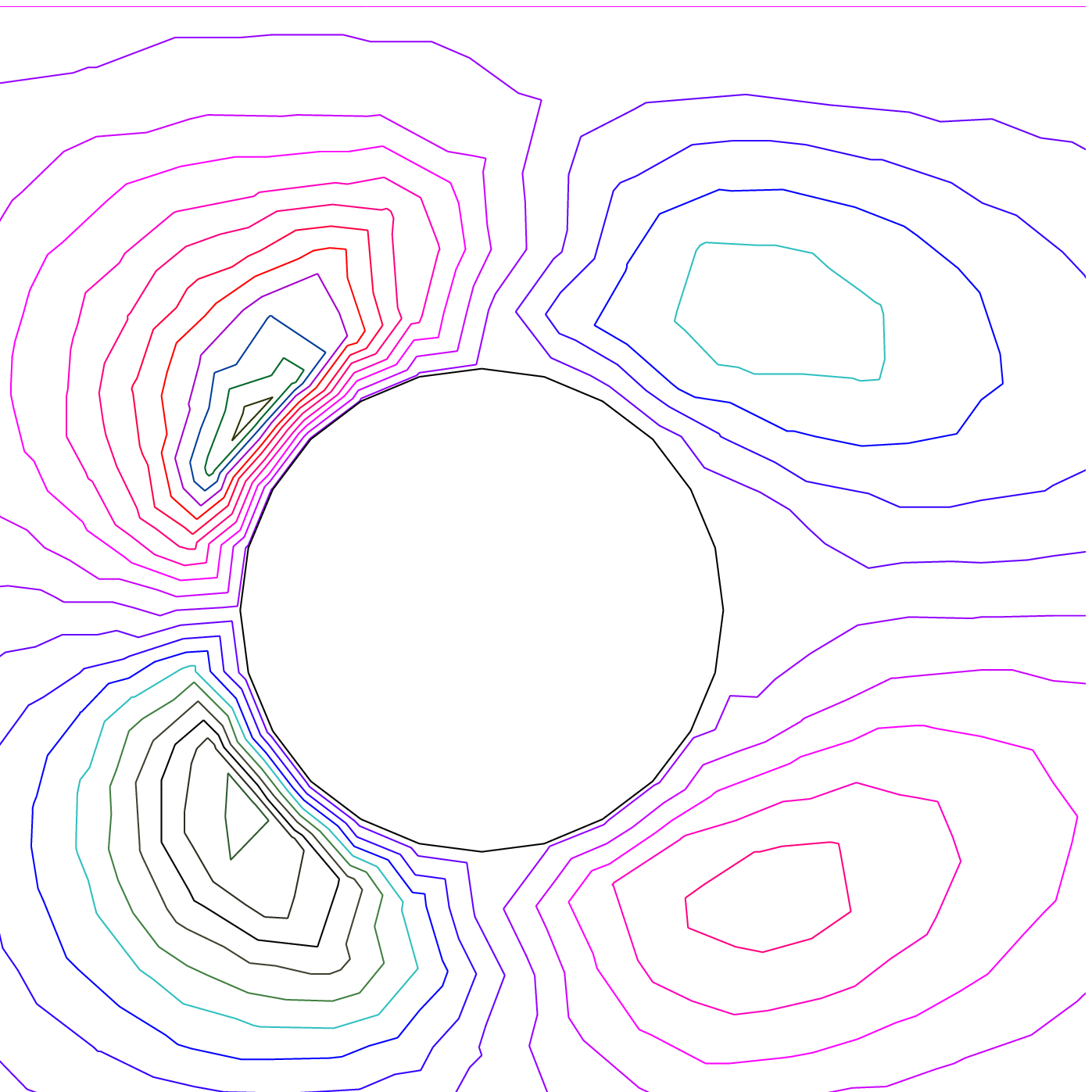}}\\
  {(b) $y_2$ for initial shape.}
  \end{minipage}
\begin{minipage}[b]{0.32\textwidth}
  \centering
 {\includegraphics[width=1.6in]{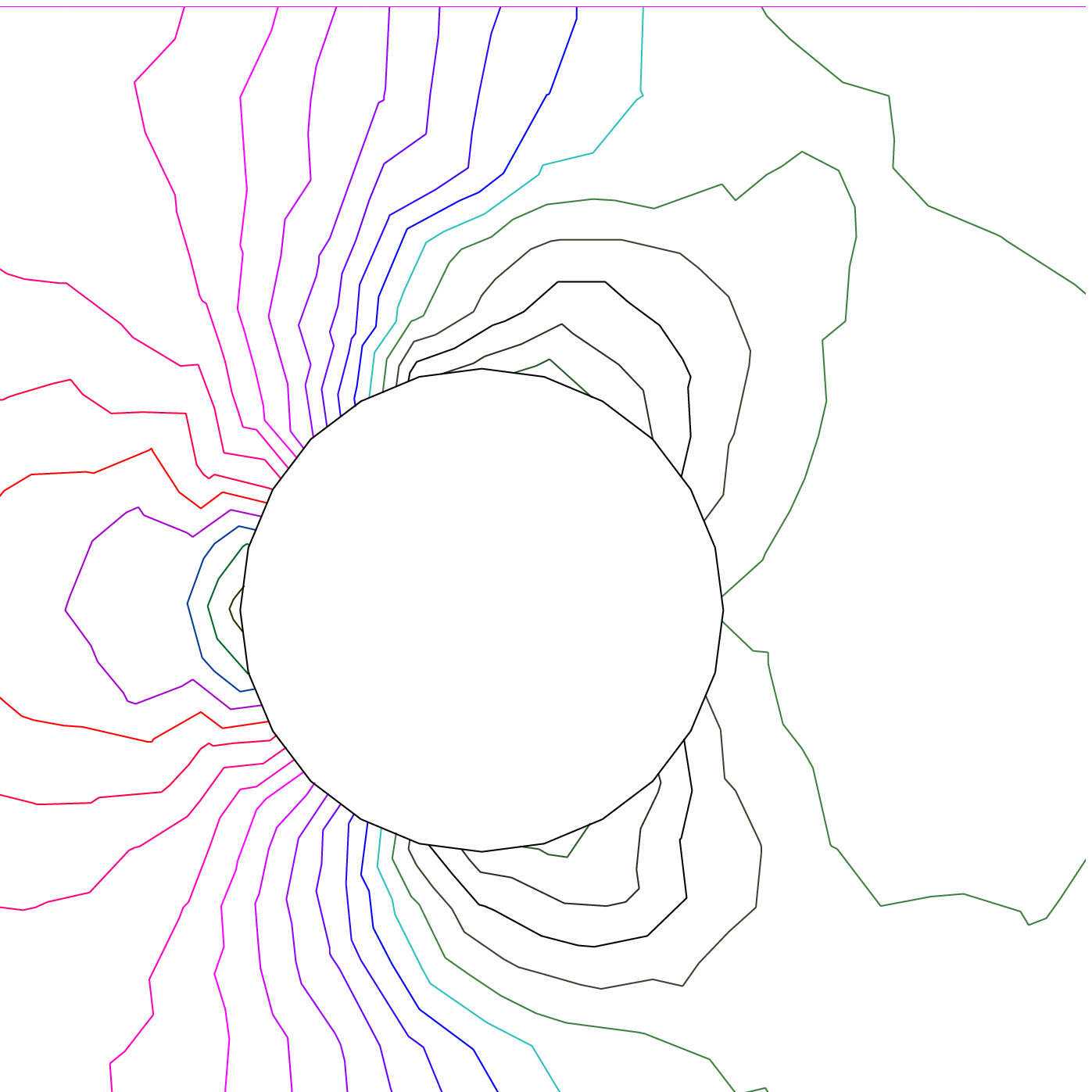}}\\
  {(c) $p$ for initial shape.}
  \end{minipage}
  \\
\begin{minipage}[b]{0.32\textwidth}
  \centering
 {\includegraphics[width=1.6in]{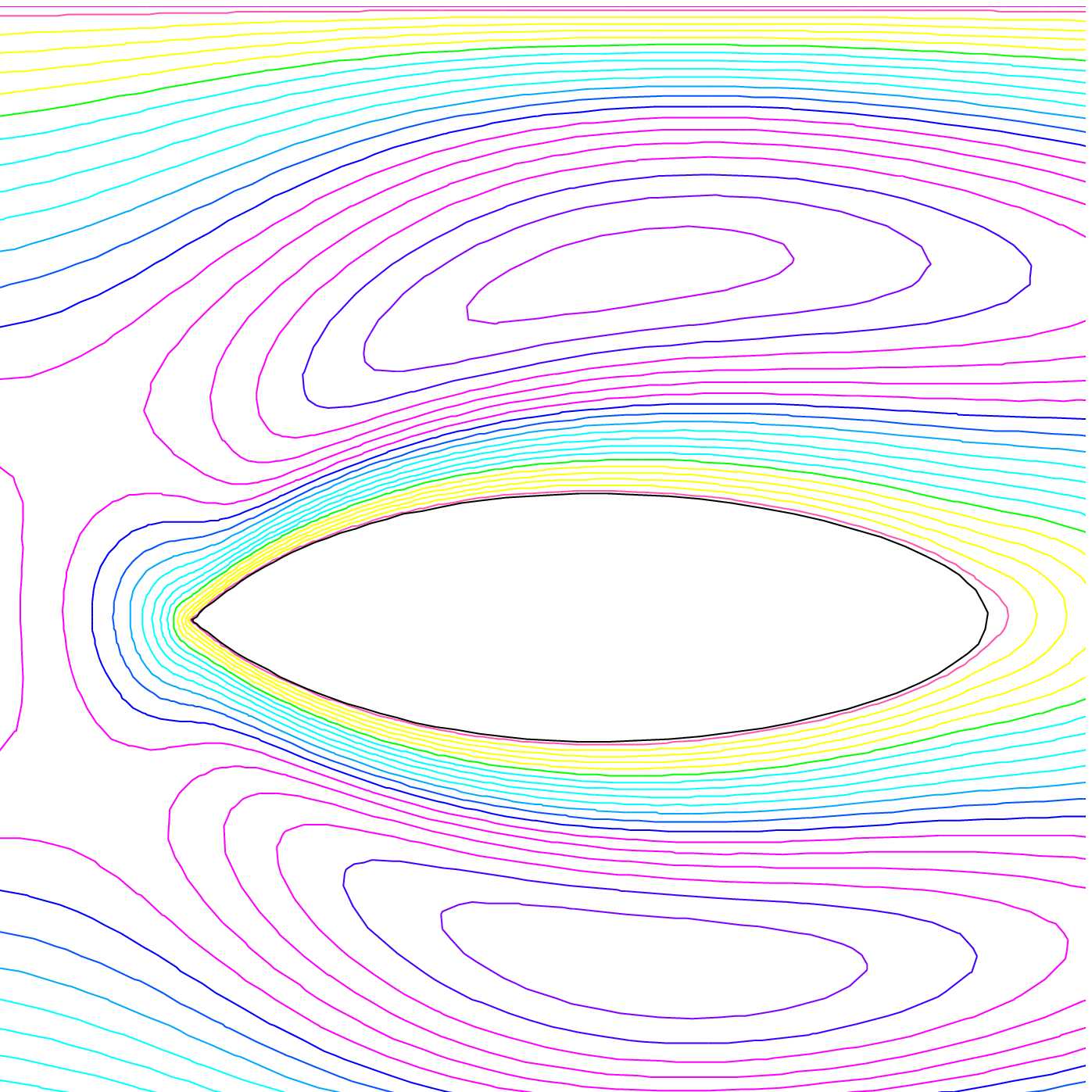}}\\
  {(d) $y_1$ for optimal shape.}
  \end{minipage}
\begin{minipage}[b]{0.32\textwidth}
  \centering
{\includegraphics[width=1.6in]{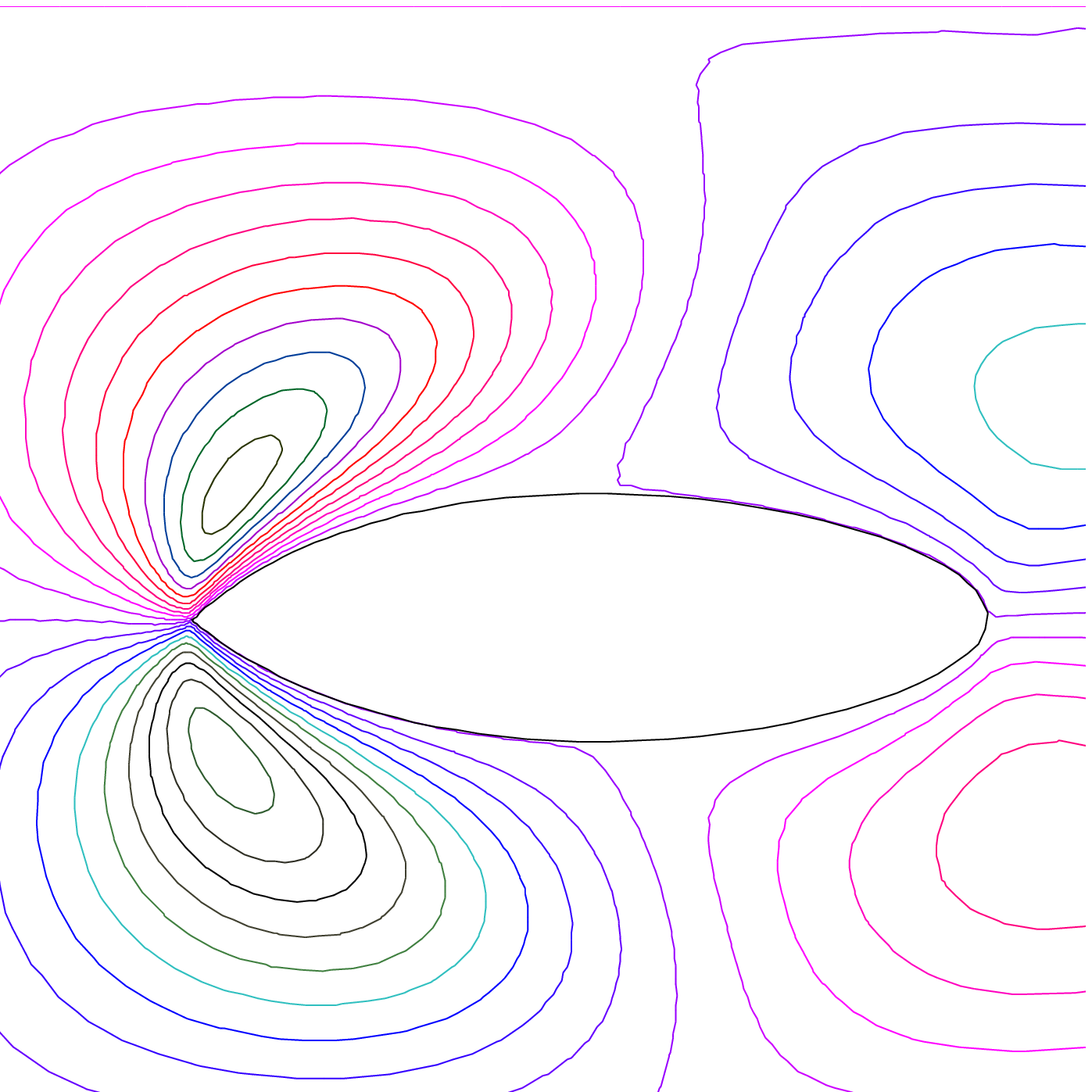}}\\
  {(e) $y_2$ for optimal shape.}
  \end{minipage}
\begin{minipage}[b]{0.32\textwidth}
  \centering
 {\includegraphics[width=1.6in]{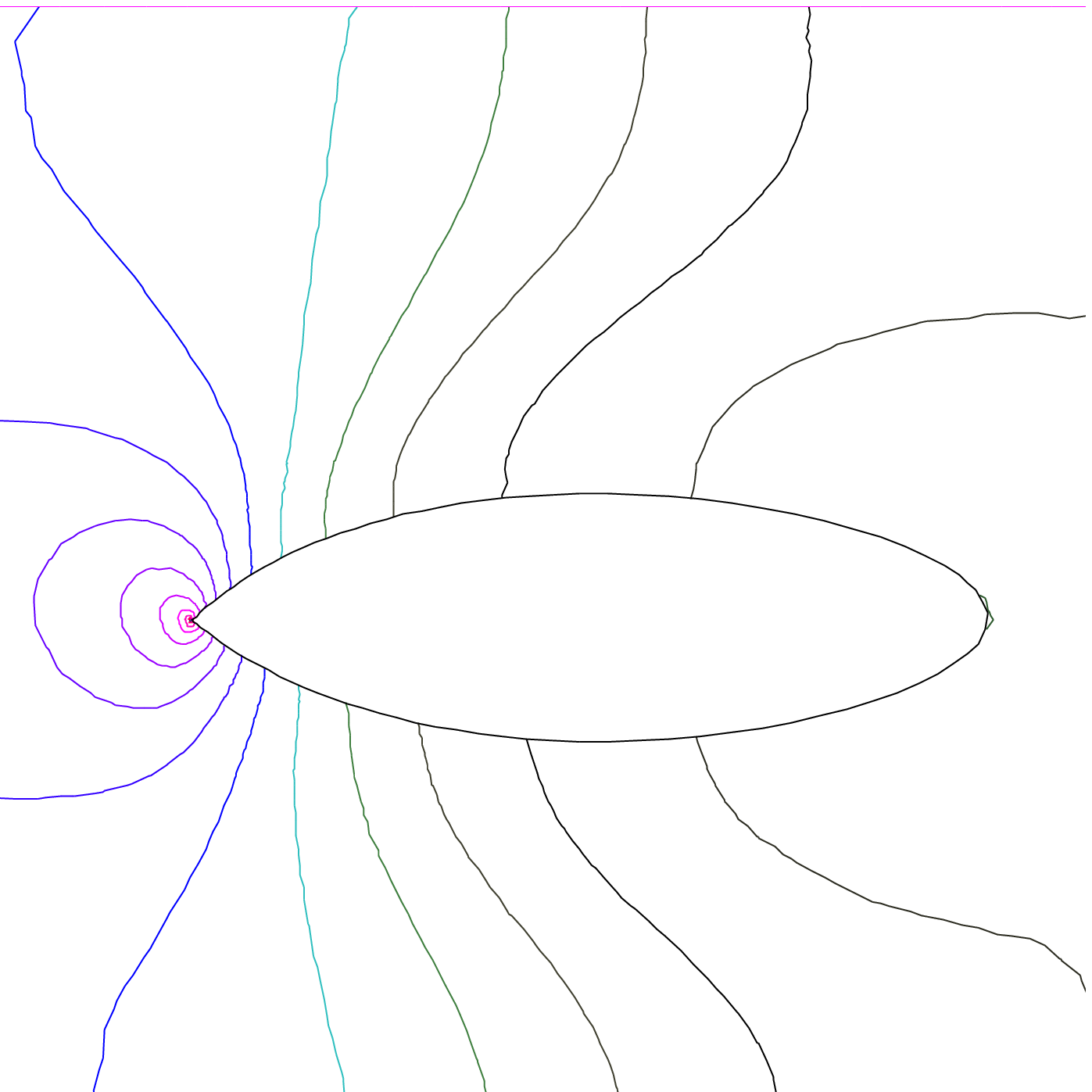}}\\
  {(f) $p$ for optimal shape.}
   \end{minipage}
 \caption{\small Comparison of the initial shape and optimal shape for $\mathrm{Re}=40.$\label{fig0:e}}
  \end{figure}
\begin{figure}[!htbp]
\renewcommand{\captionlabelfont}{\small}
\setcaptionwidth{5.5in}
\begin{minipage}[b]{0.32\textwidth}
  \centering
 {\includegraphics[width=1.6in]{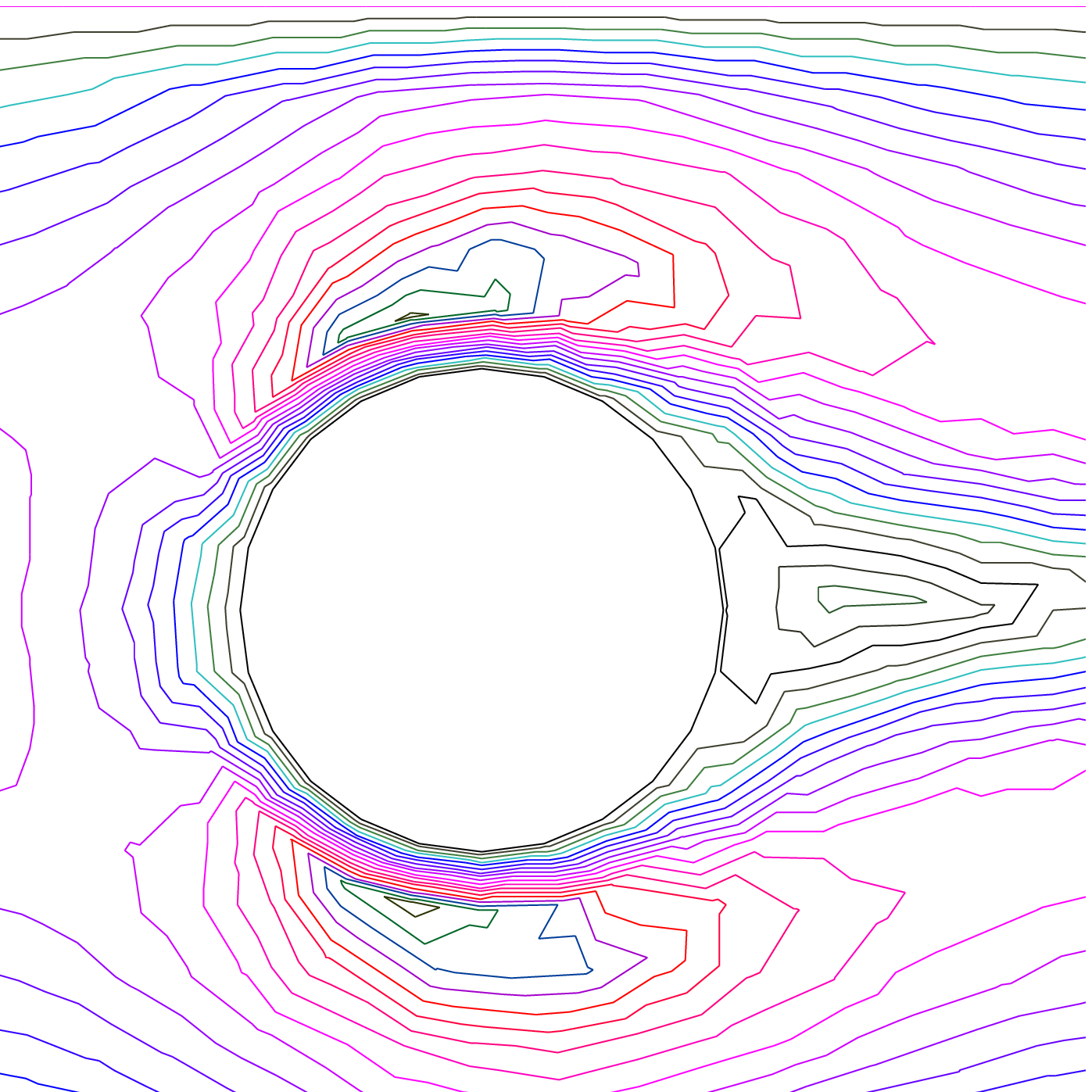}}\\
  {(a) $y_1$ for initial shape.}
  \end{minipage}
\begin{minipage}[b]{0.32\textwidth}
  \centering
{\includegraphics[width=1.6in]{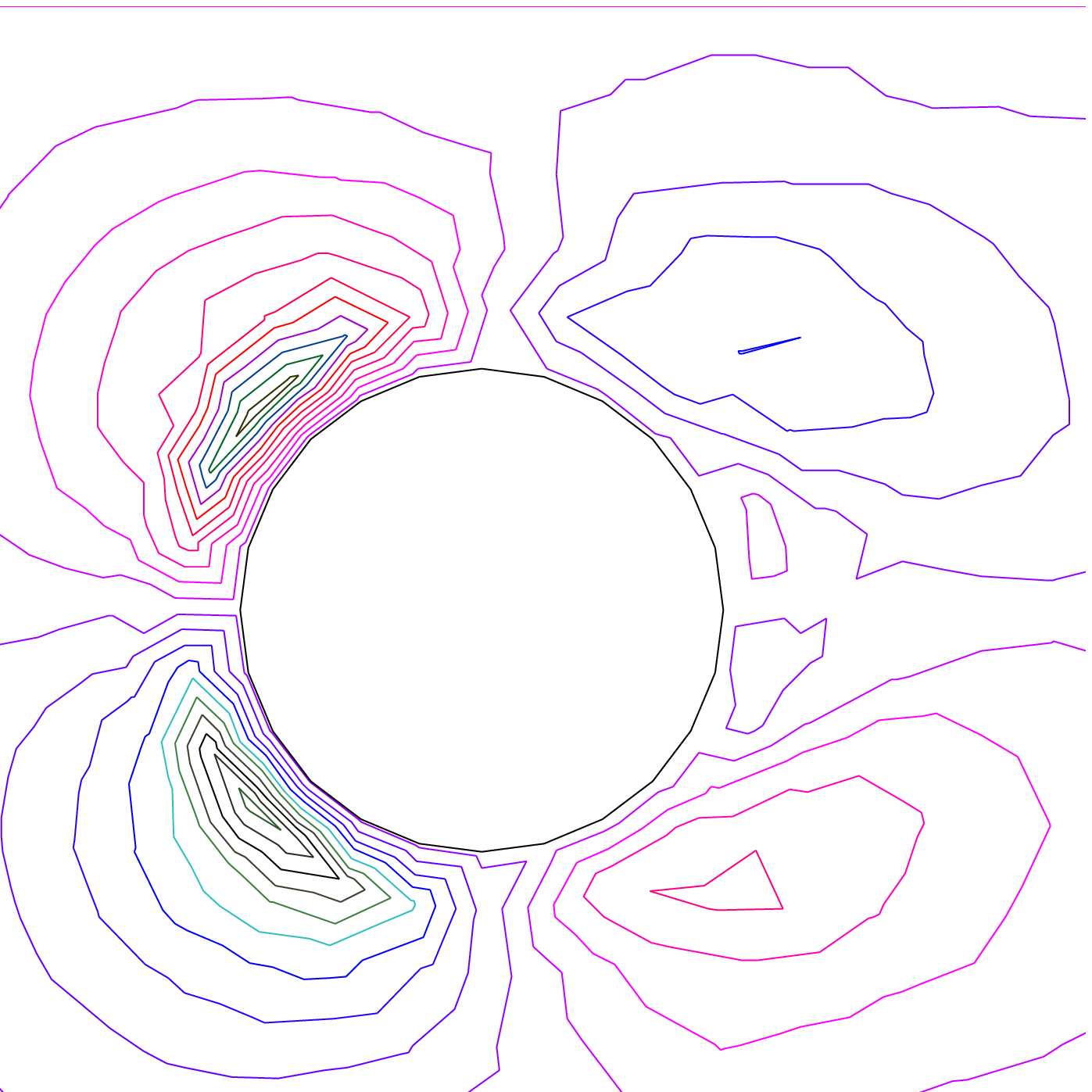}}\\
  {(b) $y_2$ for initial shape.}
  \end{minipage}
\begin{minipage}[b]{0.32\textwidth}
  \centering
 {\includegraphics[width=1.6in]{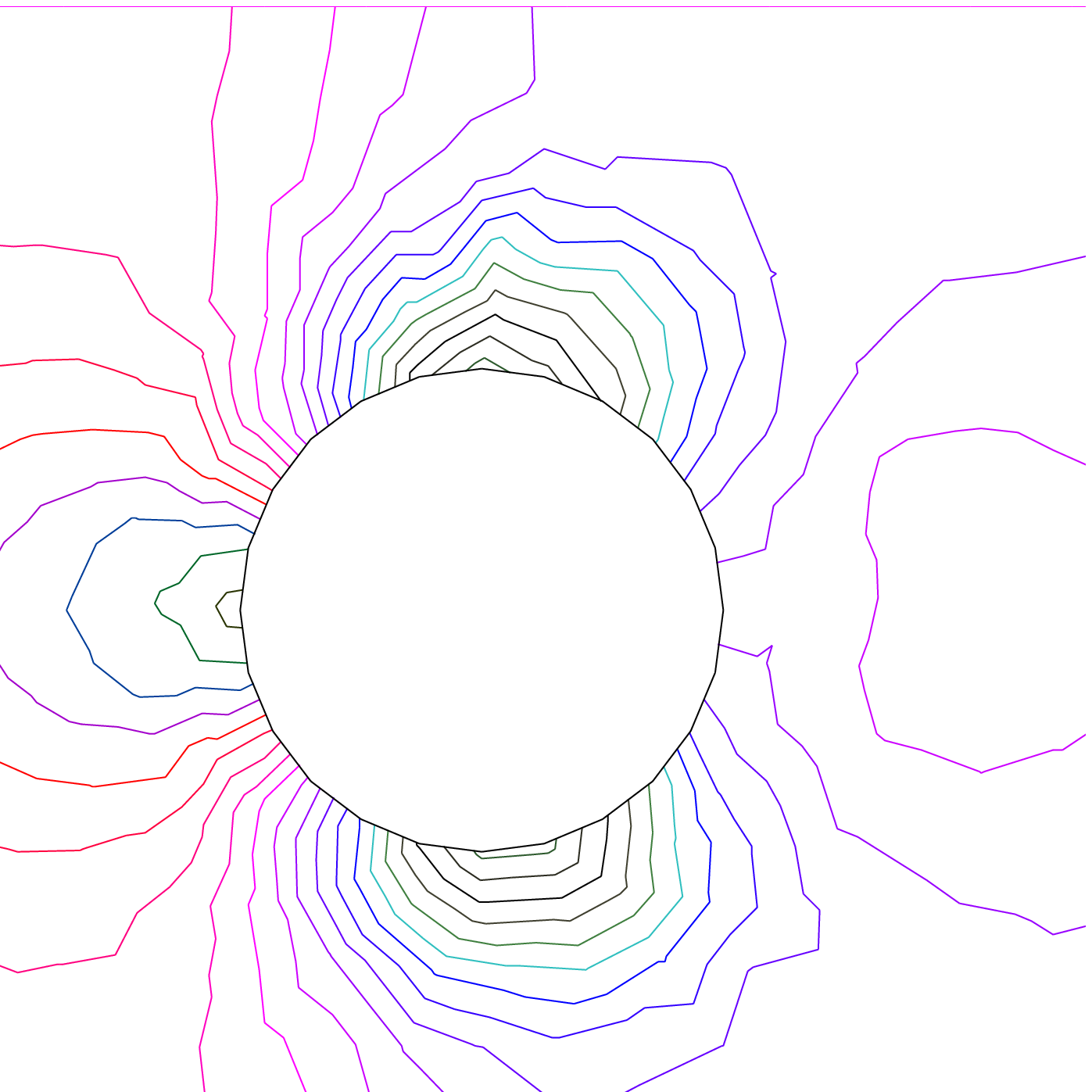}}\\
  {(c) $p$ for initial shape.}
  \end{minipage}
  \\
\begin{minipage}[b]{0.32\textwidth}
  \centering
 {\includegraphics[width=1.6in]{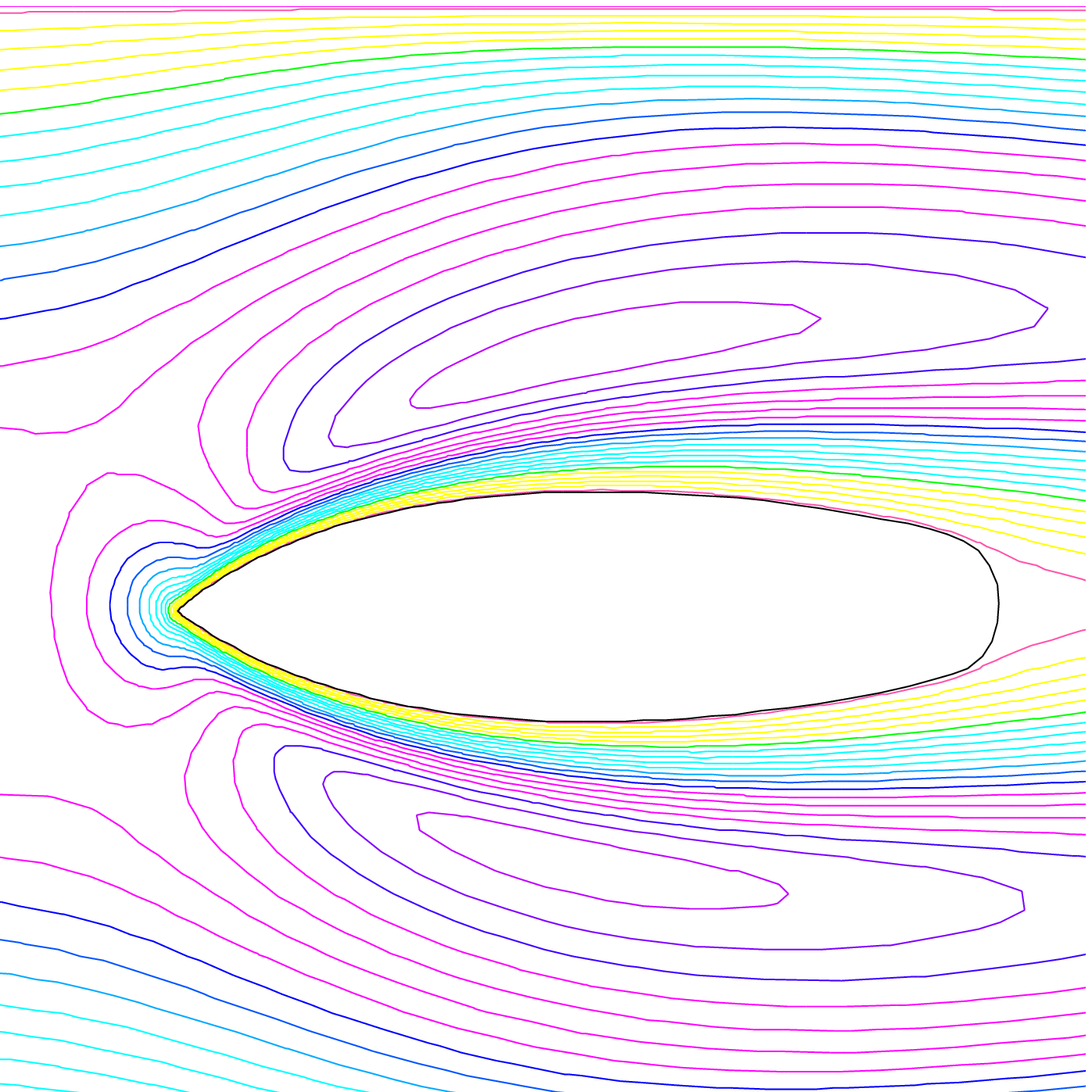}}\\
  {(d) $y_1$ for optimal shape.}
  \end{minipage}
\begin{minipage}[b]{0.32\textwidth}
  \centering
{\includegraphics[width=1.6in]{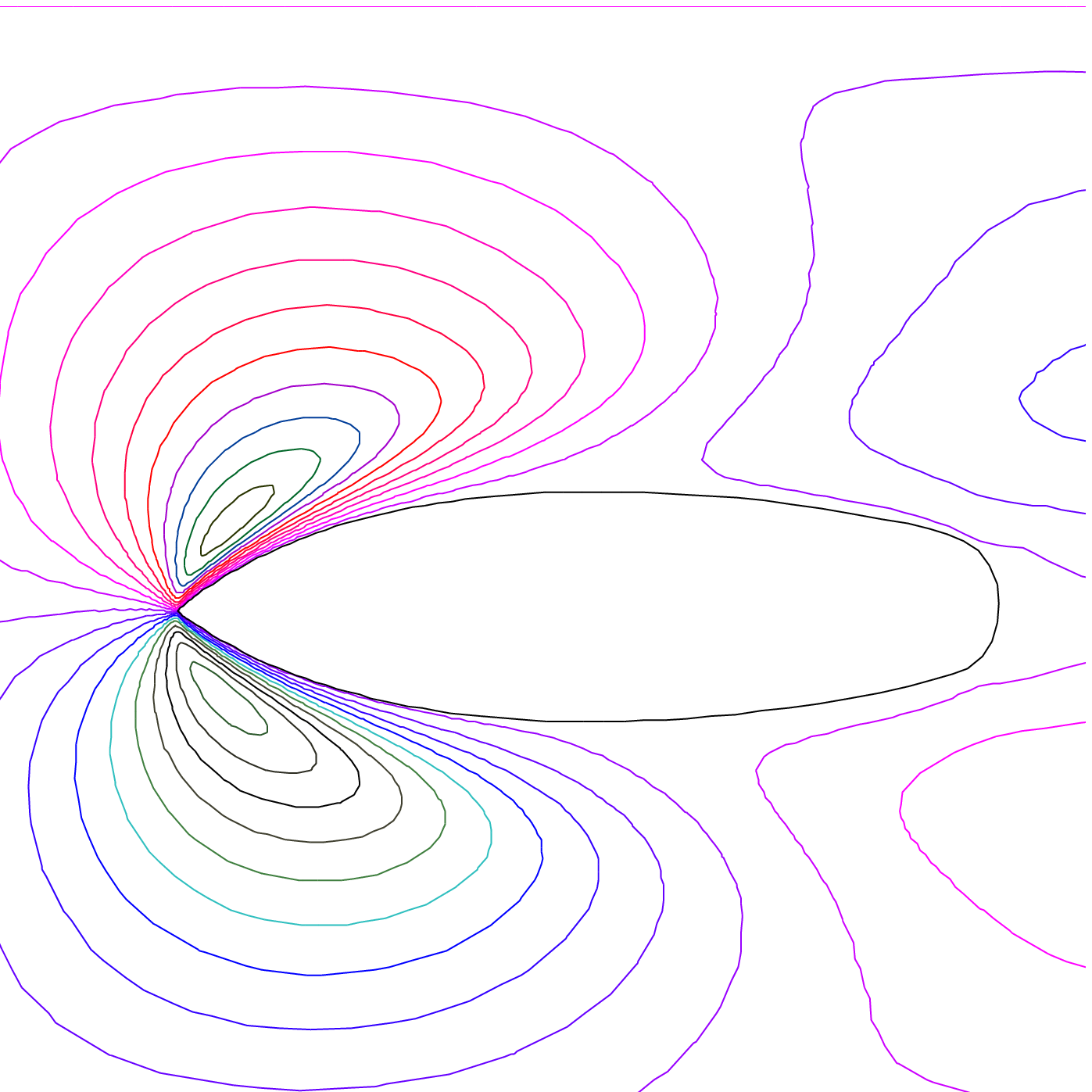}}\\
  {(e) $y_2$ for optimal shape.}
  \end{minipage}
\begin{minipage}[b]{0.32\textwidth}
  \centering
 {\includegraphics[width=1.6in]{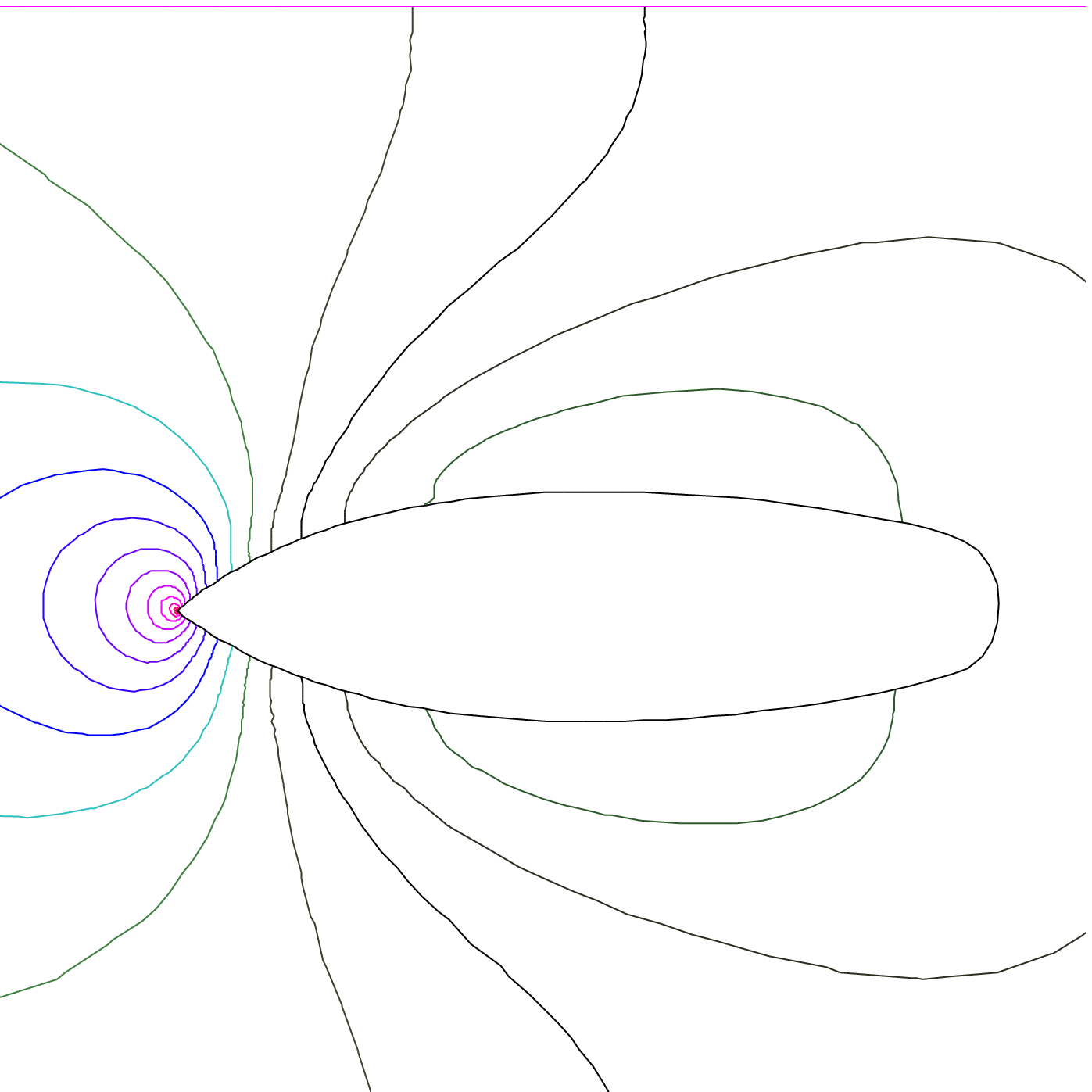}}\\
  {(f) $p$ for optimal shape.}
   \end{minipage}
 \caption{\small Comparison of the initial shape and optimal shape for $\mathrm{Re}=200$.\label{fig0:c}}
  \end{figure}
\begin{figure}[!htbp]
\renewcommand{\captionlabelfont}{\small}
\setcaptionwidth{3.4in}
  \centering
  \includegraphics[width=0.8\textwidth]{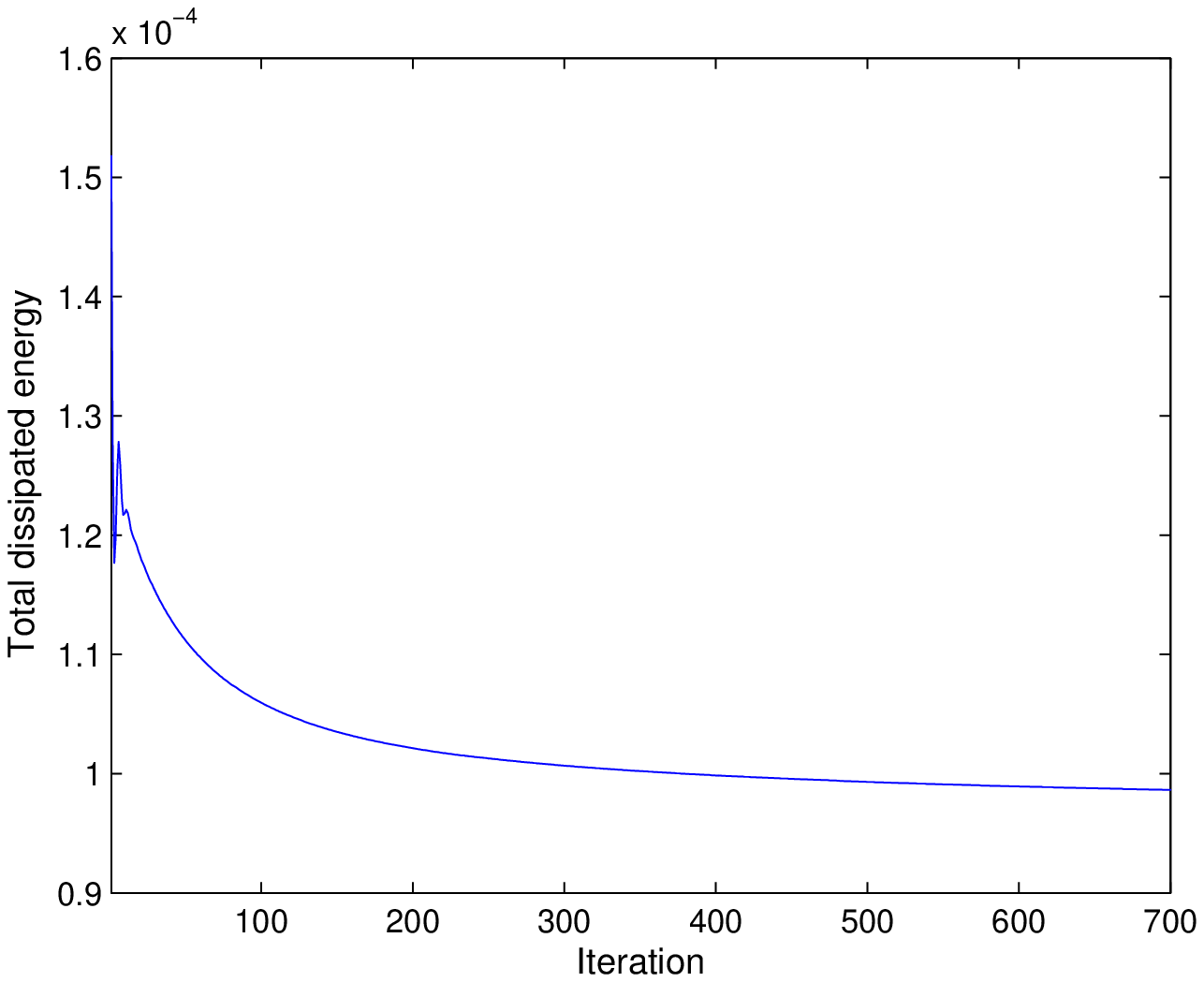}
    \caption{\small Convergence history of the total dissipated energy for $\mathrm{Re}=40$.\label{fig:cost2}}
\end{figure}
\begin{figure}[!htbp]
\renewcommand{\captionlabelfont}{\small}
\setcaptionwidth{3.4in}
  \centering
  \includegraphics[width=0.8\textwidth]{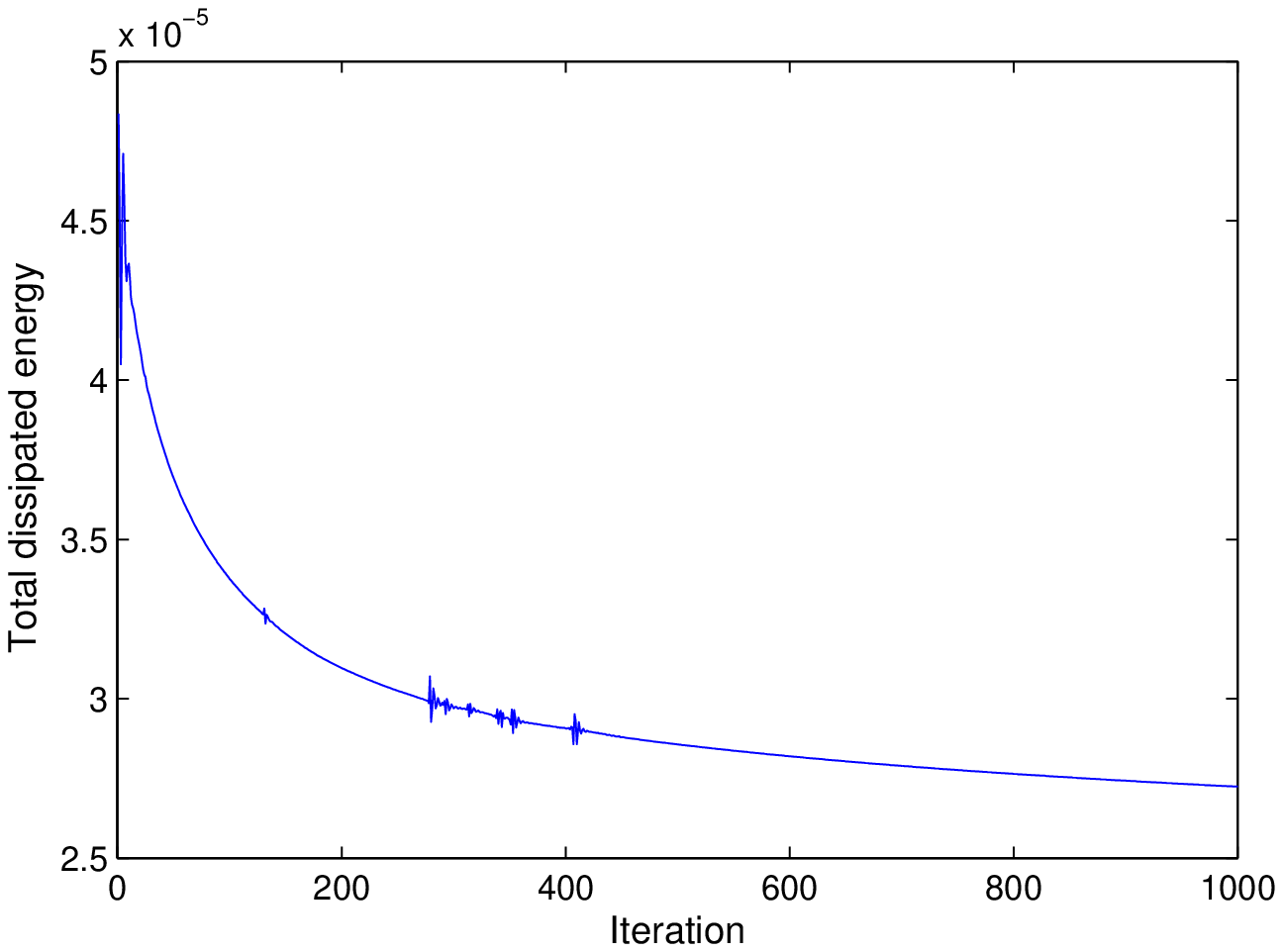}
  \caption{\small Convergence history of the total dissipated energy for
  $\mathrm{Re}=200$.\label{fig:cost3}}
\end{figure}

\end{document}